\documentclass{amsart}

\usepackage{amsfonts}
\usepackage{amsmath}
\usepackage{amsthm}
\usepackage{amssymb}
\usepackage[english]{babel}
\usepackage{graphics}
\usepackage{latexsym}
\usepackage{longtable} 
\usepackage{mathrsfs} 
\usepackage{graphicx}
\usepackage{wrapfig} %\usepackage{txfonts}
\usepackage[pdftex,colorlinks=true,linkcolor=blue,citecolor=magenta]{hyperref}
\usepackage{enumitem}
\usepackage{esint}

\newtheorem{theorem}{Theorem}
\newtheorem{corollary}[theorem]{Corollary}
\newtheorem{lemma}[theorem]{Lemma}
\newtheorem{proposition}[theorem]{Proposition}

\newtheorem{claim}[theorem]{Claim}
\newtheorem{example}[theorem]{Example}

\theoremstyle{definition}
\newtheorem{definition}[theorem]{Definition}

\newcommand{\mL}{\mathcal{L}}
\newcommand{\mH}{\mathcal{H}}
\newcommand{\mF}{\mathcal{F}}

\newcommand{\mW}{\mathcal{W}}
\newcommand{\mX}{\mathfrak{X}}

\newcommand{\mM}{\mathcal{M}}

\newcommand{\J}{\mathcal{J}}

\newcommand{\I}{\mathrm{I}}
\newcommand{\D}{\mathrm{D}}
\newcommand{\M}{\mathrm{M}}

\newcommand{\A}{\mathrm{A}}
\newcommand{\B}{\mathrm{B}}
\newcommand{\K}{\mathrm{K}}
\renewcommand{\L}{\mathrm{L}}
\renewcommand{\J}{\mathrm{J}}

\newcommand{\R}{\mathbb{R}}

\newcommand{\N}{\mathbb{N}}

\newcommand{\mB}{\mathbb{B}}

\renewcommand{\L}{\mathrm{L}}

\newcommand{\noi}{\noindent}
\newcommand{\ms}{\medskip}

\newcommand{\al}{\alpha}

\newcommand{\ga}{\gamma}

\newcommand{\e}{\varepsilon}
\newcommand{\si}{\sigma}

\newcommand{\la}{\lambda}

\newcommand{\ka}{\kappa}
\newcommand{\Om}{\Omega}
\newcommand{\om}{\omega}

\newcommand{\av}{-\hspace{-10.5pt}\displaystyle\int}

\newcommand{\weak }{\, -\!\!\!\!-\!\!\!\!\rightharpoonup}
\newcommand{\weakstar }{ \overset{\, *_{\phantom{|}}}{{\smash{\weak }}\, } }

\newcommand{\larrow}{\longrightarrow}
\newcommand{\ot}{\otimes}

\newcommand{\LL}{\text{\LARGE$\llcorner$}}
\newcommand{\p}{\partial}
\newcommand{\sub}{\subseteq}

\newcommand{\by}{\times}

\newcommand{\dist}{\mathrm{dist}}

\renewcommand{\div}{\mathrm{div}}

\newcommand{\bt}{\begin{theorem}}\newcommand{\et}{\end{theorem}}
\newcommand{\bd}{\begin{definition}}\newcommand{\ed}{\end{definition}}
\newcommand{\bl}{\begin{lemma}}\newcommand{\el}{\end{lemma}}
\newcommand{\beq}{\begin{equation}}\newcommand{\eeq}{\end{equation}}
\newcommand{\bc}{\begin{claim}}\newcommand{\ec}{\end{claim}}
\newcommand{\bex}{\begin{example}}\newcommand{\eex}{\end{example}}
\newcommand{\bcor}{\begin{corollary}}\newcommand{\ecor}{\end{corollary}}
\newcommand{\bp}{\begin{proof}}\newcommand{\ep}{\end{proof}}

\newcommand{\BPT}{\medskip \noindent \textbf{Proof of Theorem} }

\numberwithin{equation}{section}
\begin{document}

\title[A problem in ${\mathrm{L}}^\infty$ with PDE constraints]{A minimisation problem  in ${\mathrm{L}}^\infty$ with PDE and unilateral constraints}
 
\author{Nikos Katzourakis}

\address{Department of Mathematics and Statistics, University of Reading, Whiteknights, PO Box 220, Reading RG6 6AX, United Kingdom}

\email{n.katzourakis@reading.ac.uk}

  \thanks{\!\!\!\!\!\!\!\!\texttt{The author has been partially financially supported by the EPSRC grant EP/N017412/1}}
  
%\subjclass[2010]{35D99, 35D40, 35J47, 35J47, 35J92, 35J70, 35J99}

\date{}


\begin{abstract} We study the minimisation of a cost functional which measures the misfit on the boundary of a domain between a component of the solution to a certain parametric elliptic PDE system and a prediction of the values of this solution. We pose this problem as a PDE-constrained minimisation problem for a supremal cost functional in ${\mathrm{L}}^\infty$, where except for the PDE constraint there is also a unilateral constraint on the parameter. We utilise approximation by PDE-constrained minimisation problems in ${\mathrm{L}}^p$ as $p\to\infty$ and the generalised Kuhn-Tucker theory to derive the relevant variational inequalities in ${\mathrm{L}}^p$ and ${\mathrm{L}}^\infty$. These results are motivated by the mathematical modelling of the novel bio-medical imaging method of Fluorescent Optical Tomography.
\end{abstract}

\maketitle

%\tableofcontents

\section{Introduction}   \label{Section1}

Let $\Om \sub \R^n$ be an open bounded set with $ {\mathrm{C}}^1$ boundary $\p\Om$ and let also $n\geq 3$. Consider the next Robin boundary value problem for a pair of coupled linear elliptic systems:
\beq \label{1.3}
\left\{\ \ 
\begin{array}{lll}
(a)\ \ & \ \  -\div(\D u \hspace{1pt} \A)\,+\, \K u\, =\, S, & \ \ \text{ in }\Om, \ms
\\
(b)&\ \     -\div( \D v \hspace{1pt} \B)\,+\, \L v\, =\, \xi \M u, & \ \ \text{ in }\Om,  \ms
\\
(c)&  \ \ \ \ \, (\D u \hspace{1pt} \A) \hspace{1pt} \mathrm n \,+\, \ga u\, =\, s, & \ \ \text{ on }\p\Om,  \ms
\\
(d)& \ \  \ \ \, (\D v \hspace{1pt} \B) \hspace{1pt} \mathrm n \,+\, \ga v\, =\, 0, & \ \ \text{ on }\p\Om,
\end{array}
\right.
\eeq
where $u,v : \Om \larrow \R^2$ are the solutions, $\mathrm n : \p\Om \larrow \R^n$ is the outer unit normal vector field on $\p\Om$ and the coefficients $\A,\B,\K,\L,\M,s,S,\xi,\ga$ satisfy $\ga>0$ and
\beq 
\label{1.6}
\left\{
\begin{split}
u,\, v,\, S\ & : \ \Om\larrow \R^2, \ \ \ \ \ \D u,\, \D v\ : \ \Om \larrow \R^{2\by n},
\\
\K,\,\L,\,\M\ & : \ \Om \larrow \R^{2\by 2}, \ \ \ \ \ \, \A,\,\B\ :\ \Om \larrow \R^{n\by n}_+,
\\
s\ & : \ \p \Om\larrow \R^2, \ \ \ \ \ \ \ \ \ \ \ \  \xi\  : \ \Om\larrow [0,\infty).
\end{split}
\right.
\eeq
Here the matrix-valued maps $\K,\L$ are assumed to have the form
\beq  \label{1.5}
\K\,:=\,
\left[
\begin{array}{cc}
k_1 & -k_2
\\
k_2 & k_1
\end{array}
\right]
,\ \ \
\L\,:=\,
\left[
\begin{array}{cc}
l_1 & -l_2
\\
l_2 & l_1
\end{array}
\right].
\eeq
We will suppose that there exists $a_0>0$ such that
\beq 
\label{3.1}
\left\{
\begin{split}
&\A,\B \in \mathrm{VMO}(\R^n;\R^{n\by n}_+),\ \ \si(\A),\si(\B) \sub \Big[a_0,\frac{1}{a_0} \Big],
\\
& \K,\L,\M \in {\mathrm{L}}^{\infty}(\Om;\R^{2\by 2}),\ \ \ \, k_1,\, \geq\, a_0,\ \ \ l_1\, \geq\, a_0.
\end{split}
\right.
\eeq
We note that our general notation will be either standard or self-explanatory, as e.g.\ in the textbooks \cite{E,KV}. The PDE problem above (where the coefficient $\xi$ is considered as a parameter) arises in Fluorescent Optical Tomography, which is a new and evolving bio-medical imaging method with wider industrial applications. Optical tomography is being very intensely studied, as it presents some pros over standard imaging methods which use X-rays, Gamma-rays, electromagnetic radiation and ultrasounds. In particular, it is less harmful for living organisms and more precise. In this setting, the coefficients take the form
\[
\begin{split}
\A(x) \, =\, \B(x) \, &=\, \frac{1}{3\big(\mu_{ai}(x) \,+\, \mu'_{s}(x)  \,+\, \xi(x) \big)} \mathrm I_3,
\\
\K(x)\,=\,\L(x)\, &=\, \left[
\begin{array}{cc}
\mu_{ai}(x)\,+\, \xi(x) & -\dfrac{\om}{\mathrm c}
\\
\dfrac{\om}{\mathrm c} & \mu_{ai}(x)\,+\, \xi(x)
\end{array}
\right],
\\
\M(x)\,&=\, \left[
\begin{array}{cc}
\dfrac{\phi}{1+(\om \tau)^2} & -\dfrac{\phi\om \tau}{1+(\om \tau)^2}
\\
\dfrac{\phi\om \tau}{1+(\om \tau)^2} & \dfrac{\phi}{1+(\om \tau)^2}
\end{array}
\right],
\end{split}
\]
where $\mathrm I_3$ is the identity matrix in $\R^3$. In the above, the coefficients $\A,\B$ describe the diffusion of photons, $\mu'_s$ is the reduced scattering coefficient, $\phi$ is the quantum efficiency of the fluorophore, $\mu_{ai}$ is the absorption coefficient due to the endogenous chromophores, $\tau$ is the fluorophore lifetime, $\xi$ is the absorption coefficient due to the exogenous fluorophore and $\om$ is the modulated light frequency and $\mathrm c$ the speed of light. Finally, $S,s$ are light sources. 

Technically, the aim of optical tomography is to reconstruct $\xi$ in the domain from measurements of light intensity on the boundary. A fluorescent dye is injected into $\Om$. In order to determine the dye concentration $\xi$, $\Om$ is illuminated by a light source $s=s(x)$ placed on the boundary $\p\Om$. The light diffuses inside $\Om$, and wherever dye is present, infrared light is emitted that can then be detected again on the boundary through appropriate sensors. The goal is then to reconstruct $\xi$ from the obtained boundary images. For more details we refer to \cite{ARH, AMM, Ar, BJ, BJ2, FP, FES, GHA, GEZ, JBS, ZRWN, ZG}.

In this paper, motivated by the problem of optical tomography and by the recent developments in Calculus of Variations in $L^\infty$ appearing in the papers \cite{KM,KP1,KPa}, we consider the problem of minimising over the class of all admissible parameters $\xi$ a certain cost functional which measures the deviation of the solution $v$ on the boundary $\p\Om$ from some prediction $\tilde v$ of its values. Given the high complexity of the optical tomography problem, in this work which is the companion paper of \cite{K5} we will make the simplifying assumption that the diffusion coefficients $\A,\B$ and the optical terms $\K,\L$ do not depend explicitly on the dye distribution. On the other hand, though, we allow for potentially different diffusions and coefficients in the two systems describing the red and the infrared light. To this end, fix $N\in \N$, $m \in [n,\infty)$ and $p >\max\left\{n,{2n}/({n-2}) \right\}$. Consider Borel sets
\beq \label{4.1}
\big\{\mathrm{B}_1,...,\mathrm{B}_N\big\}  \sub \p\Om
\eeq
and maps 
\beq  \label{4.2}
\big\{s_1,...,s_N\big\} \sub {\mathrm{L}}^{\frac{m}{2}}(\p\Om;\R^{2}),  \ \ \ \ \big\{S_1,...,S_N \big\} \sub {\mathrm L}^{\frac{nm}{2n+m}}(\Om;\R^{2}),
\eeq
and let also 
\beq
 \label{4.3}
\big\{\tilde v_1,...,\tilde v_N\big\}  \sub {\mathrm{L}}^{\infty}(\p\Om;\R^{2})
\eeq
be predicted (noisy) values of the solution $v$ of \eqref{1.3}(b)-\eqref{1.3}(d) on the boundary $\p\Om$. Suppose that for any $i\in\{1,...,N\}$, the pair $(u_i,v_i)$ solves \eqref{1.3} with coefficiens $(S_i,s_i,\xi)$. For the $N$-tuple of solutions $(u_1,...,u_N;v_1,...,v_N)$, we will symbolise
\[
\big(\vec u,\vec v \, \big) \, \in \, {\mathrm{W}}^{1,\frac{m}{2}}(\Om;\R^{2\by N}) \by {\mathrm{W}}^{1,p}(\Om;\R^{2 \by N})
\]
and understand $(u_i)_{i=1...N}$ and $(u_i)_{i=1...N}$ as matrix valued. Similarly, we will see the corresponding vectors of test functions as
\[
\big(\vec \phi,\vec \psi \, \big) \, \in \, {\mathrm{W}}^{1,\frac{m}{m-2}}(\Om;\R^{2\by N}) \by {\mathrm{W}}^{1,{\frac{p}{p-1}}}(\Om;\R^{2 \by N}).
\]
Our aim is to determine some $\xi \in {\mathrm{L}}^p(\Om,[0,\infty))$ such that all the misfits
\[
\big|(v_i - \tilde v_i)\big|_{\mathrm{B}_i}\big|
\]
between the predicted approximate solution and the actual solution are minimal. We will minimise the error in ${\mathrm{L}}^\infty$ by means of approximations in ${\mathrm{L}}^p$ for large $p$ and then take the limit $p\to \infty$. By minimising in ${\mathrm{L}}^\infty$ one can achieve uniformly small cost, rather than on average. Since no reasonable cost functional is coercive in our admissible class, we will therefore follow two different approaches to rectify this problem, but in a unified fashion. The first and more popular idea is to add a Tykhonov-type regularisation term $ \al \| \xi \|$ for small $\al>0$ and some appropriate norm. The alternative approach is to consider that an a priori ${\mathrm{L}}^\infty$ bound is given on $\xi$. The latter approach appears to be more natural for applications, as it does not alter the error functional. For finite $p<\infty$, we can relax this to an ${\mathrm{L}}^p$ bound, but as we are mostly interested in the limit case $p=\infty$, we will only discuss the case of ${\mathrm{L}}^\infty$ bound. In view of the above observations, we define the integral functional
\beq 
\label{4.4}
\I_p\big(\vec u,\vec v ,\xi \big) \,:=\, \sum_{i=1}^N \big\|v_i - \tilde v_i\big\|_{\dot{\mathrm{L}}^p(\mathrm{B}_i)}  +\, \al \| \xi \|_{\dot{\mathrm{L}}^p(\Om)}, \ \ \ \ (\vec u,\vec v \,,\xi \big) \in \mathfrak{X}^{p}(\Om)
\eeq
and its supremal counterpart
\beq
\label{4.5}
\I_\infty\big(\vec u,\vec v ,\xi \big) \,:=\, \sum_{i=1}^N \big\|v_i - \tilde v_i\big\|_{{\mathrm{L}}^\infty(\mathrm{B}_i)}  +\, \al \| \xi \|_{{\mathrm{L}}^\infty(\Om)} \ \ \ \ (\vec u,\vec v \,,\xi \big) \in \mathfrak{X}^{\infty}(\Om),
\eeq
where the dotted $\dot{\mathrm{L}}^p$ quantities are regularisations of the respective norms:
\beq
\label{4.6}
\|f \|_{\dot{\mathrm{L}}^p(\Om)} := \left(\, \av_{\Om}(|f|_{(p)})^p\, \mathrm d \mL^{n}\! \right)^{\!\!1/p}, \ \ \  \|g \|_{\dot{\mathrm{L}}^p(\mathrm{B}_i)} := \left(\, \av_{\mathrm{B}_i}(|g|_{(p)})^p\, \mathrm d \mH^{n-1} \!\right)^{\!\!1/p}.
\eeq
The slashed integrals denote the average with respect to the Lebesgue measure $\mL^n$ and the Hausdorff measure $\mH^{n-1}$ respectively and  $|\cdot |_{(p)}$ is a regularisation of the Euclidean norm away from zero in the corresponding space, given by
\beq
\label{4.7}
| \cdot |_{(p)} \,:=\, \sqrt{| \cdot |^2 + p^{-2}}.
\eeq
The admissible classes $\mX^p(\Om)$ and $\mX^\infty(\Om)$ into which we will minimise \eqref{4.4}-\eqref{4.5} are defined by setting
\beq
\label{4.10}
\phantom{\Big|}\mathscr{X}^p(\Om)\,:=\, {\mathrm{W}}^{1,\frac{m}{2}}(\Om;\R^{2\by N}) \by {\mathrm{W}}^{1,p}(\Om;\R^{2 \by N}) \by {\mathrm{L}}^p(\Om),
\eeq
\beq
\label{4.8}
\mX^p(\Om)\,:=\, \left\{ 
\begin{array}{l}
 (\vec u,\vec v \,,\xi \big) \in \mathscr{X}^p(\Om) :  \text{ for all } i \in\{1,...,N\},\ (u_i,v_i,\xi) \text{ satisfies }\ms\ms
\\
\hspace{80pt} 0\leq \xi \leq M \, \text{ a.e. on $\Om$}
\\
\text{and}\ms
\\
\ \  \left\{\ \ 
\begin{array}{lll}
(a)_i\ \ & \ \  -\div(\D u_i  \hspace{1pt} \A)\,+\, \K u_i\, =\, S_i, & \ \ \text{ in }\Om, \ms
\\
(b)_i&\ \     -\div( \D v_i  \hspace{1pt} \B)\,+\, \L v_i\, =\, \xi \M u_i, & \ \ \text{ in }\Om,  \ms
\\
(c)_i&  \ \ \ \ (\D u_i  \hspace{1pt} \A)\hspace{1pt} \mathrm n\,+\, \ga u_i\, =\, s_i, & \ \ \text{ on }\p\Om,  \ms
\\
(d)_i& \ \ \ \  ( \D v_i  \hspace{1pt} \B)\hspace{1pt} \mathrm n\,+\, \ga v_i\, =\, 0, & \ \ \text{ on }\p\Om,
\end{array}
\right.
\\
\\
 \text{for }\A,\B,\K,\L,\M,S_i,s_i,\xi,\ga,p\text{ satisfying hypotheses \eqref{1.6}-\eqref{4.3}}
\end{array}\!\!
\right\}
\eeq
and
\beq
\label{4.9}
\mX^\infty(\Om)\, :=\, \bigcap_{n<p<\infty} \mX^p(\Om).
\eeq
Note that $\mX^\infty(\Om)$ is a subset of a Frech\'et space, rather than of a Banach space, but no difficulties will emerge out of this. We will assume that the pair of parameters $\al,M \in [0,\infty]$ satisfy {\it either of the two conditions:}
\begin{align}
a>0 \ \ \ \text{ and }\ \  M=\infty \label{4.11},
\\
a=0 \ \ \ \text{ and }\ \  M<\infty \label{4.12}.
\end{align}
Namely, if we have non-zero regularisation parameter $\al$, then no upper ${\mathrm{L}}^\infty$ bound $M$ is required, but if the parameter $\al$ vanishes we need an upper ${\mathrm{L}}^\infty$ bound $M$.

Our first main result concerns the existence of $\I_p$-minimisers in $\mX^p(\Om)$, the existence of $\I_\infty$-minimisers in $\mX^\infty(\Om)$ and the approximability of the latter by the former as $p\to \infty$.

\begin{theorem}[$\I_\infty$-misfit minimisers, $\I_p$-misfit minimisers \& convergence as $p\to\infty$] \label{th4}

\ms

\noi {\rm (A)} The functional $\I_p$ has a constrained minimiser $(\vec u_p,\vec v_p ,\xi_p)$ in the admissible class $\mX^p(\Om)$:
\beq \label{4.15A}
\I_p\big(\vec u_p,\vec v_p ,\xi_p \big)\, =\, \inf\Big\{\I_p\big(\vec u,\vec v ,\xi \big)\, : \ \big(\vec u,\vec v ,\xi \big) \in \mX^p(\Om) \Big\}.
\eeq
{\rm (B)} The functional $\I_\infty$ has a constrained minimiser $(\vec u_\infty,\vec v_\infty ,\xi_\infty)$ in the admissible class $\mX^\infty(\Om)$
\beq  \label{4.15B}
\I_\infty\big(\vec u_\infty,\vec v_\infty ,\xi_\infty \big)\, =\, \inf\Big\{\I_\infty\big(\vec u,\vec v ,\xi \big)\, : \ \big(\vec u,\vec v ,\xi \big) \in \mX^\infty(\Om) \Big\}.
\eeq
Additionally, there exists a subsequence of indices $(p_j)_1^\infty$ such that the sequence of respective $\I_{p_j}$-minimisers $\big(\vec u_{p_j},\vec v_{p_j} ,\xi_{p_j} \big)$ satisfy as $p_j\to \infty$ that
\beq
\label{4.16}
\left\{ \ \
\begin{array}{ll}
\xi_p \weak \xi_\infty, & \!\!\!\!\!\!\!\!\!\!\!\! \text{in }{\mathrm{L}}^q(\Om), \text{ for all }q\in (1,\infty),
\ms
\\
\vec u_{p} \weak \vec u_\infty, & \!\!\!\!\!\!\!\!\!\!\!\! \text{in }{\mathrm{W}}^{1,\frac{m}{2}}(\Om;\R^{2\by N}), \ms
\\
\vec u_{p} \larrow \vec u_\infty, & \!\!\!\!\!\!\!\!\!\!\!\! \text{in }{\mathrm{L}}^{\frac{m}{2}}(\Om;\R^{2\by N}), \ms
\\
\vec v_{p} \weak \vec v_\infty, & \!\!\!\!\!\!\!\!\!\!\!\!\text{in }{\mathrm{W}}^{1,q}(\Om;\R^{2\by N}), \text{ for all }q\in (1,\infty),\ms
\\
\vec v_p \larrow \vec v_\infty, & \!\!\!\!\!\!\!\!\!\!\!\! \text{in }  {\mathrm{C}}^0(\overline{\Om};\R^{2\by N}),  
\ms
\\
\I_p\big(\vec u_{p},\vec v_{p} ,\xi_{p} \big) \larrow & \!\!\! \I_\infty\big(\vec u_\infty,\vec v_\infty ,\xi_\infty\big). 
\end{array}
\right.
\eeq
\end{theorem}

Given the existence of constrained minimisers established by Theorem \ref{th4} above, the next natural question concerns the existence of necessary conditions in the form of PDEs governing the constrained minimisers. Unlike the case of unconstrained minimisation, no analogue of Euler-Lagrange equations is available in this case. One the one hand, the PDE constraints will give rise to Lagrange multipliers which are functionals. On the other hand, the unilateral constraint on $\xi$ gives rise to a variational differential inequality, rather than an equation. This generalised variational context of extrema with constraints is known as the Kuhn-Tucker theory (see e.g.\ \cite{Z}). Hence, our next main result regarding the variational inequalities for finite $p$ is given below.

\begin{theorem}[Variational inequalities in ${\mathrm{L}}^p$] \label{th7} For any $p>\max\{n,2n/(n-2)\}$, there exist Lagrange multipliers
\[
\big(\vec \phi_p, \vec \psi_p \big) \ \in {\mathrm{W}}^{1,\frac{m}{m-2}} (\Om;\R^{2\by N}) \by {\mathrm{W}}^{1,\frac{p}{p-1}} (\Om;\R^{2\by N})
\]
associated with the constrained minimisation problem \eqref{4.15A}, such that $\big(\vec u_p, \vec v_p, \xi_p\big) \in \mX^p(\Om)$ satisfies the relations
\beq
\label{5.20}
\begin{split}
& \int_\Om (\eta-\xi_p) \left( \al\, \frac{\mathrm d [\mu_p(\xi_p)]}{\mathrm{d}\mL^n} \, + \sum_{i=1}^N  \big(\M  u_{pi}  \big)  \cdot \psi_{pi} \right) \, \mathrm{d}\mL^n \, \geq\, 0,
\end{split}
\eeq
\beq
\label{5.21}
\begin{split}
 \int_{\p\Om} \vec w : \mathrm d [\vec \nu_p(\vec v_p)] \, =\,\sum_{i=1}^N & \bigg\{  \int_\Om \Big[ \B :(\D w_i^\top \D \psi_{pi}) \, +  \big(\L w_i\big) \cdot \psi_{pi} \Big] \, \mathrm{d}\mL^n 
\\
& \ + \int_{\p\Om} (\ga w_i)\cdot \psi_{pi} \, \mathrm{d}\mH^{n-1} \bigg\},
\end{split}
\eeq
\beq
\label{5.22}
\begin{split}
  \sum_{i=1}^N & \bigg\{\int_\Om\Big[ \A :(\D z_i^\top \D \phi_{pi}) \, + (\K z_i)\cdot \phi_{pi} \Big] \, \mathrm{d}\mL^n \,+ \int_{\p\Om}(\ga z_i)\cdot \phi_{pi} \, \mathrm{d}\mH^{n-1} \bigg\}
\\
& \ \ =\ \sum_{i=1}^N \int_\Om   \xi_p \big(\M z_i \big) \cdot \psi_{pi} \, \mathrm{d}\mL^n  , 
\end{split}
\eeq
for any test functions $\big(\vec v, \vec w, \eta\big) \in \mathscr{X}^p(\Om)$. In \eqref{5.20}-\eqref{5.21}, $\mu_p(\xi)$ is the next $\xi$-dependent real Radon measure in $\mM(\Om;\R)$
\beq 
\label{5.3}
\mu_p(\xi)\, :=\, \frac{ (|\xi|_{(p)})^{p-2} \xi}{ \mL^{n}(\Om) \big(\|\xi \|_{\dot{\mathrm{L}}^p(\Om)}\big)^{p-1}} \mL^{n}\LL_{\Om} ,
\eeq
and $\vec \nu_p(\vec v)$ is the next $\vec v$-dependent matrix-valued Radon measure in $\mM\big(\p\Om;\R^{2\by N}\big)$ 
\beq 
\label{5.2}
\vec \nu_p(\vec v)\, :=\, \sum_{i=1}^N \Bigg(\frac{\big(|v_i-\tilde v_i|_{(p)}\big)^{p-2} (v_i-\tilde v_i)}{ \mH^{n-1}(\mathrm{B}_i) \big(\big\|v_i-\tilde v_i\big\|_{\dot{\mathrm{L}}^p(\mathrm{B}_i)}\big)^{p-1}}  \ot e_i \Bigg)  \mH^{n-1}\LL_{\mathrm{B}_i} .
\eeq
\end{theorem}

Note that the measures $\vec \nu_p(\vec v)$ and $\mu_p(\xi)$ are absolutely continuous with respect to the Hausdorff measure $\mH^{n-1}\LL_{\p\Om}$ and the Lebesgue measure $\mL^n\LL_{\Om}$ respectively. In general, $\mu_p(\xi)$ is signed if $\xi \in {\mathrm{L}}^p(\Om)$, but due to the constraint we have $\mu_p(\xi_p)\geq 0$. Further,  $\{e_1,...,e_N\}$ symbolises the standard Euclidean basis of $\R^N$, ``$:$" symbolises the standard inner product in $\R^{2\by N}$, and ${\mathrm d [\mu_p(\xi_p)]}/{\mathrm{d}\mL^n}$ symbolises the Radon-Nikodyn derivative of $\mu_p(\xi_p)$ with respect to $\mL^n$. The reason we obtain three different relations of which one is inequality and two are equations is the following. If we ignore the PDE constraints in \eqref{4.8}, then the admissible class is the Cartesian product of two vector spaces (spaces for $\vec u$ and $\vec v$), and a convex set (space of $\xi$, see also \eqref{5.10} that follows). Since the unilateral constraint is only for $\xi$, the variational inequality arises only for this variable. The decoupling of these relations is a consequence of linear independence.
 
Our final main result is the limiting counterpart of Theorem \ref{th7} for the constrained minimiser of the $\infty$-problem. To this aim, let us set
\[
C_\infty:=\, \limsup_{p_j \to \infty}\, C_p,\ \ \ \ C_p\,:=\, \|\vec \phi_p\|_{{\mathrm{W}}^{1,\frac{m}{m-2}} (\Om)} +\,  \| \vec \psi_p \|_{{\mathrm{W}}^{1,1}(\Om)},
\]
where $\big(\vec \phi_p, \vec \psi_p \big)$ are the Lagrange multipliers associated with the constrained minimisation problem \eqref{4.15B} (Theorem \ref{th7}).

\begin{theorem}[Variational inequalities in ${\mathrm{L}}^\infty$] \label{th14}

If additionally $m>n$, there exists a subsequence $(p_j)_1^\infty$ and a pair of limiting measures 
\[
(\mu_\infty,\vec \nu_\infty) \, \in \mM\big(\Om;[0,\infty)\big) \by \mM\big(\p\Om;\R^{2\by N}\big)
\]
such that
\beq
\label{6.2}
\big(\mu_p(\xi_p) , \vec \nu_p(\vec v_p)\big) \, \weakstar \ \big(\mu_\infty, \vec \nu_\infty\big) \ \text{ in } \ \mM(\Om;\R) \by \mM\big(\p\Om;\R^{2\by N}\big),
\eeq
as $p_j \to \infty$. Then:

\ms

\noi {\rm (I)} If $C_\infty=0$ and $\al>0$, then $\xi_\infty =0$ a.e.\ on $\Om$ \ and \ $\vec \nu_\infty = \vec 0$.

\ms

\noi {\rm (II)} If $C_\infty>0$, then there exist (rescaled) limiting Lagrange multipliers
\[
\big(\vec \phi_\infty, \vec \psi_\infty \big) \ \in {\mathrm{W}}^{1,\frac{m}{m-2}} (\Om;\R^{2\by N}) \by \mathrm{BV} \big(\Om;\R^{2\by N}\big)
\]
such that
\beq
\label{6.1}
 \big( C_p^{-1}\vec \phi_p \, , \, C_p^{-1}\vec \psi_p\big) \ \weakstar \ \big(\vec \phi_\infty, \vec \psi_\infty \big)
\eeq
as $p_j \to \infty$, in the space ${\mathrm{W}}^{1,\frac{m}{m-2}} (\Om;\R^{2\by N}) \by \mathrm{BV} (\Om;\R^{2\by N})$. In this case, the constrained minimiser $\big(\vec u_\infty, \vec v_\infty, \xi_\infty\big) \in \mX^\infty(\Om)$ satisfies the next three relations:
\beq
\label{6.3}
 \frac{\al}{C_\infty} \int_\Om \eta\, \mathrm d \mu_\infty \, 
+\  \sum_{i=1}^N \int_\Om (\eta-\xi_\infty) \big(\M  u_{\infty i}  \big)  \cdot \psi_{\infty i} \, \mathrm{d}\mL^n \, \geq\, \frac{\al}{C_\infty}\| \xi_\infty\|_{{\mathrm{L}}^\infty(\Om)},
\eeq
\beq
\label{6.4}
\begin{split}
 \frac{1}{C_\infty}\int_{\p\Om} \vec w : \mathrm d \vec \nu_\infty \, =\,\sum_{i=1}^N & \bigg\{  \int_\Om  \B : (\D w_i)^\top  \mathrm{d}[\D \psi_{\infty i}] \, 
 +  \int_{\p\Om}  \big(\L w_i\big) \cdot \psi_{\infty i}  \, \mathrm{d}\mL^n 
\\
& \ + \int_{\p\Om} (\ga w_i)\cdot \psi_{\infty i} \, \mathrm{d}\mH^{n-1} \bigg\},
\end{split}
\eeq
\beq
\label{6.5}
\begin{split}
  \sum_{i=1}^N & \bigg\{\int_\Om\Big[ \A :(\D z_i^\top \D \phi_{\infty i}) \, + (\K z_i)\cdot \phi_{\infty i} \Big] \, \mathrm{d}\mL^n \,+ \int_{\p\Om}(\ga z_i)\cdot \phi_{\infty i} \, \mathrm{d}\mH^{n-1} \bigg\}
\\
& \ \ =\ \sum_{i=1}^N \int_\Om   \xi_\infty \big(\M z_i \big) \cdot \psi_{\infty i} \, \mathrm{d}\mL^n  , 
\end{split}
\eeq
for any 
\[
\big(\vec z, \vec w,\eta\big) \, \in \, {\mathrm{C}}^1(\overline{\Om};\R^{2\by N}) \by {\mathrm{C}}_0^1(\overline{\Om};\R^{2\by N})\by {\mathrm{C}}^0(\overline{\Om};[0,M]). 
\]

\end{theorem}

We conclude this lengthy introduction with some comments about the general variational context we use herein. Calculus of Variations in ${\mathrm{L}}^\infty$ is a modern subarea of analysis pioneered by Aronsson in the 1960s (see \cite{A1}-\cite{A4}) who considered variational problems of supremal functionals, rather than integral functional. For a pedagogical introduction we refer e.g.\ to \cite{C,K4}. Except for their endogenous mathematical appeal, ${\mathrm{L}}^\infty$ cost-error functionals are important for applications because by minimising their supremum rather than their average (as e.g.\ in standard ${\mathrm{L}}^2$ approaches), we obtain improved performance/predictions/fitting. Indeed, minimisation of the supremum of the misfit guarantees uniform smallness, namely  deviation spikes of small volume are excluded. Interesting results regarding ${\mathrm{L}}^\infty$ variational problems can be found e.g.\ in \cite{AB, BBJ, BaJe, BJW1, BN, BP, CDP, GNP, MWZ, PP, P, RZ}.

\ms

\section{Proofs}  

We begin with an auxiliary result of independent interest, namely the well-posedness of general Robin boundary value problems for linear systems. 

\begin{theorem}[Well-posedness in ${\mathrm{W}}^{1,2}$ and ${\mathrm{W}}^{1,p}$] \label{th1} Let $\Om \Subset \R^n$ be a domain with $ {\mathrm{C}}^1$ boundary and let $\mathrm n : \p\Om \larrow \R^n$ be the outer unit normal. Consider the boundary value problem 
\beq \label{2.1}
\left\{\ \ 
\begin{array}{ll}
  -\div(\D u \hspace{1pt} \A)\,+\, \K u\, =\, f-\div F, & \ \ \text{ in }\Om, \ms
\\
\ \ (\D u \hspace{1pt} \A - F)\hspace{1pt} \mathrm n\,+\, \ga u\, =\, g, & \ \ \text{ on }\p\Om,
\end{array}
\right.
\eeq
where  $\ga>0$. We suppose there exists $a_0>0$ such that
\beq  \label{2.2}
\left\{
\ \ 
\begin{split}
& \A \in {\mathrm{L}}^\infty(\Om;\R^{n\by n}_+),\ \ \ a_0 |z|^2 \leq  \A:z \ot z  \leq \frac{1}{a_0}|z|^2\ \ \forall \, z\in\R^n ,
\\
& \K \in {\mathrm{L}}^\infty(\Om;\R^{2\by 2}),\ \ \ \ \K\,:=\,
\left[
\begin{array}{cc}
k_1 & -k_2
\\
k_2 & k_1
\end{array}
\right]
\  \text{ and } \ \ k_1\geq a_0.
\end{split}
\right.
\eeq
If
\beq \label{2.3}
f \in {\mathrm{L}}^2(\Om;\R^2),\ \ \ F \in {\mathrm{L}}^2(\Om;\R^{2\by n}),\ \ \ g \in {\mathrm{L}}^2(\p\Om;\R^2),
\eeq
then, \eqref{2.1} has a unique weak solution in ${\mathrm{W}}^{1,2}(\Om;\R^2)$ satisfying
\beq \label{2.4}
\left\{\ \ 
\begin{split}
\int_\Om\Big[ \A :(\D u^\top \D\phi ) \, +\,  (\K u)\cdot \phi  \Big] \, \mathrm{d}\mL^n \,+\, \int_{\p\Om}\big[ \ga u\cdot\phi  \big] \, \mathrm{d}\mH^{n-1}
\\
= \int_\Om\Big[ f\cdot\phi \, +\,  F:\D \phi \Big] \, \mathrm{d}\mL^n\, +\, \int_{\p\Om}\big[ g \cdot \phi \big] \, \mathrm{d}\mH^{n-1},
\end{split}
\right.
\eeq   
for all $\phi \in {\mathrm{W}}^{1,2}(\Om;\R^2)$. In addition, exists $C>0$ depending only on the coefficients such that
\beq 
\label{2.5}
\| u \|_{{\mathrm{W}}^{1,2}(\Om)} \, \leq\, C\Big(\| f \|_{{\mathrm{L}}^2(\Om)} +\, \| F \|_{{\mathrm{L}}^2(\Om)} +\,\| g \|_{{\mathrm{L}}^{2}(\p\Om)}\Big).
\eeq
If additionally for some $p > {2n}/({n-2})$ we have
\[
\A \in \mathrm{VMO}(\R^n;\R^{n\by n}_+),\ \ f \in {\mathrm{L}}^{\frac{np}{n+p}}(\Om;\R^2),\ \ F \in {\mathrm{L}}^p(\Om;\R^{2\by n}),\ \ g \in {\mathrm{L}}^p(\p\Om;\R^2),
\]
then, the weak solution of \eqref{2.1} lies in the space ${\mathrm{W}}^{1,p}(\Om;\R^2)$. In addition, there exists $C>0$ depending only on the coefficients and $p$ such that
\beq \label{2.7}
\| u \|_{{\mathrm{W}}^{1,p}(\Om)} \, \leq\, C\Big(\| f \|_{{\mathrm{L}}^{\frac{np}{n+p}}(\Om)} +\, \| F \|_{{\mathrm{L}}^p(\Om)} +\,\| g \|_{{\mathrm{L}}^p(\p\Om)}\Big).
\eeq
\end{theorem}

In the proofs that follow  we will employ the standard practice of denoting by $C$ a generic constant whose value might change from step to step in an estimate.

\bp The aim is to apply of the Lax Milgram theorem. (Note that the matrix $\K$ is not symmetric, thus this is not a direct consequence of the Riesz theorem.) We define the bilinear functional
\[
\mB \ :\ \ {\mathrm{W}}^{1,2}(\Om;\R^2) \by {\mathrm{W}}^{1,2}(\Om;\R^2) \larrow \R,
\]
\[
\mB[u,\psi] \,:=\, \int_\Om\Big[ \A :(\D u^\top \D \psi) \, +\,  (\K u)\cdot \psi \Big] \, \mathrm{d}\mL^n \,+\, \int_{\p\Om}\big[ \ga u\cdot \psi\big] \, \mathrm{d}\mH^{n-1}.
\]
Since $\A,\K$ are ${\mathrm{L}}^\infty$, by H\"older inequality we immediately have
\[
\big| \mB[u,\psi]  \big|\,\leq\, C \| u \|_{{\mathrm{W}}^{1,2}(\Om)} \| \psi \|_{{\mathrm{W}}^{1,2}(\Om)} 
\]
for some $C>0$ and all $u,\psi \in {\mathrm{W}}^{1,2}(\Om;\R^2)$. Further, since
\[
\begin{split}
(\K u)\cdot u\, &=\, [u_1,\, u_2] \left[
\begin{array}{cc}
k_1 & -k_2
\\
k_2 & k_1
\end{array}
\right] \left[
\begin{array}{c}
u_1
\\
u_2
\end{array}
\right] \, =\, k_1 |u|^2 \, \geq \, a_0 |u|^2,
\end{split}
\]
we estimate
\[
\begin{split}
\mB[u,u] \,& \geq\, a_0\Big(\| \D u \|^2_{{\mathrm{L}}^2(\Om)} +\, \| u \|^2_{{\mathrm{L}}^2(\Om)}\Big) +\,\ga \| u \|^2_{{\mathrm{L}}^{2}(\p\Om)} ,
\end{split}
\]
for any $u \in {\mathrm{W}}^{1,2}(\Om;\R^2)$. Hence, the bilinear form $\mB$ is continuous and coercive, thus the hypotheses of the Lax-Milgram theorem are satisfied (see e.g.\ \cite{E}). Hence, for any $\Phi \in ({\mathrm{W}}^{1,2}(\Om;\R^2))^*$, exists a unique $u \in {\mathrm{W}}^{1,2}(\Om;\R^2)$ such that
\[
\ \ \ \mB[u,\psi] = \langle \Phi, \psi \rangle, \ \ \ \text{for all}\ \psi \in {\mathrm{W}}^{1,2}(\Om;\R^2).
\]
Next, we show that the functional $\Phi$ given by
\[
\langle \Phi, \psi \rangle\, :=  \int_{\p\Om}\big[ g \cdot \psi\big] \, \mathrm{d}\mH^{n-1} \, + \int_\Om\Big[ f\cdot \psi \, +\,  F:\D \psi\Big] \, \mathrm{d}\mL^n
\]
lies in $({\mathrm{W}}^{1,2}(\Om;\R^2))^*$ and we will also establish the $\mathrm L^2$ and the $\mathrm L^p$ estimates. Indeed, by the trace theorem in ${\mathrm{W}}^{1,2}(\Om;\R^2)$, there is a $C>0$ which allows to estimate
\[
\begin{split}
\big| \langle \Phi, \psi \rangle \big| \, & \leq\, \| g \|_{{\mathrm{L}}^2(\p\Om)} \| \psi \|_{{\mathrm{L}}^2(\p\Om)} \, + \Big(\| f \|_{{\mathrm{L}}^2(\Om)} +\, \| F \|_{{\mathrm{L}}^2(\Om)} \Big) \| \psi \|_{{\mathrm{W}}^{1,2}(\Om)}  
\\
& \leq\, C\Big(\| f \|_{{\mathrm{L}}^2(\Om)} +\, \| F \|_{{\mathrm{L}}^2(\Om)} +\, \| g \|_{{\mathrm{L}}^2(\p\Om)}\Big)  \| \psi \|_{{\mathrm{W}}^{1,2}(\Om)}.
\end{split}
\]
The particular choice of $\psi:=u$ together with Young inequality yield
\[
\begin{split}
\big| \langle \Phi, u \rangle \big| \, & \leq\, \e \| u \|^2_{{\mathrm{W}}^{1,2}(\Om)} \, +\, \frac{ {\mathrm{C}}^2}{4\e}\Big(\| f \|_{{\mathrm{L}}^2(\Om)} +\, \| F \|_{{\mathrm{L}}^2(\Om)} +\, \| g \|_{{\mathrm{L}}^2(\p\Om)}\Big)^{\!2}.
\end{split}
\]
We conclude with the claimed $\mathrm L^2$ estimate by combining the above estimate with our lower bound on $\mB[u,u]$.

Now we turn to the higher integrability of the weak solution. The main ingredient is to apply a well-know estimate for the Robin boundary value problem which has the form \eqref{2.7}, but applies to the scalar version of \eqref{2.1} for $\K\equiv0$, see e.g.\ \cite{ACGG, D, DK, Ge, Gh, KLS, N, N2}. Hence, we need to show that it is still true in the general case of \eqref{2.1}. To this end, we rewrite \eqref{2.1} componentwise as
\[
\ \ \left\{ \ \ 
\begin{array}{ll}
  -\div(\D u_i \A) \, =\, \big\{ f_i -\big(\K u)_i \big\}-\div F_i, & \ \ \text{ in }\Om, \ms
\\
\ \ (\D u_i \A - F_i)^\top \mathrm n\,+\, \ga u_i\, =\, g_i, & \ \ \text{ on }\p\Om,
\end{array}
\right.
\]
for $i=1,2$. By applying the estimate to the each of the components separately, we have
\beq
\label{2.8}
\begin{split}
\| u_i \|_{{\mathrm{W}}^{1,p}(\Om)} \, \leq\, & C\Big( \|\K\|_{{\mathrm{L}}^\infty(\Om)} \| u \|_{{\mathrm{L}}^{\frac{np}{n+p}}(\Om)} +\, \| f_i \|_{{\mathrm{L}}^{\frac{np}{n+p}}(\Om)}
\\
&\ \ \ +\, \| F_i \|_{{\mathrm{L}}^p(\Om)} +\,\| g_i \|_{{\mathrm{L}}^p(\p\Om)}\Big),
\end{split}
\eeq
for $i=1,2$. Note now that since we have assumed $p >{2n}/({n-2})$, we have $2 < {np}/({n+p}) < p$. Hence, by the ${\mathrm{L}}^p$ interpolation inequalities, we can estimate
\[
\  \ \ \| u \|_{{\mathrm{L}}^{\frac{np}{n+p}}(\Om)} \leq\,  \| u \|^\la_{{\mathrm{L}}^{2}(\Om)} \,  \| u \|^{1-\la}_{{\mathrm{L}}^{p}(\Om)}, \ \ \ \text{ for }\ \la = \frac{2p}{n(p-2)}.
\]
By Young's inequality
\beq
\label{2.10}
ab\, \leq \, \left\{\frac{r-1}{r}(\e r)^{\frac{1}{1-r}}\right\}b^{\frac{r}{r-1}}  \, +\, \, \e a^r,
\eeq
which holds for $a,b,\e> 0$, $r>1$ and $r/(r-1)=r'$, the choice $r:=1/(1-\la)$ yields
\[
1-\la= \frac{n(p-2)}{p(n-2)-2n},\ \ \ r=\frac{n(p-2)}{p(n-2)-2n}\ ,\ \ \ \frac{r}{r-1}=\frac{n(p-2)}{2p},
\]
and hence we can estimate 
\begin{align}
 \| u \|_{{\mathrm{L}}^{\frac{np}{n+p}}(\Om)} \, & \leq \,    \left(\| u \|_{{\mathrm{L}}^{p}(\Om)}\right)^{\frac{p(n-2)-2n}{n(p-2)}} \left(\| u \|_{{\mathrm{L}}^{2}(\Om)}\right)^{\frac{2p}{n(p-2)}} 
\nonumber\\
& \leq\, \Big(\! \left(\| u \|_{{\mathrm{L}}^{p}(\Om)}\right)^{\frac{p(n-2)-2n}{n(p-2)}} \! \Big)^{\!r} + \left[ \frac{r-1}{r}(\e r)^{\frac{1}{1-r}}\right]  \! \Big( \!\left(\| u \|_{{\mathrm{L}}^{2}(\Om)}\right)^{\frac{2p}{n(p-2)}}\! \Big)^{\!\frac{r}{r-1}}
\label{2.11}
 \\
 & =\, \e \| u \|_{{\mathrm{L}}^{p}(\Om)} \, + \left[ \frac{2p}{n(p-2)}\left( \frac{\e n (p-2)}{p(n-2)-2n} \right)^{ \!\! -\frac{p(n-2)-2n}{2p} }\right] \| u \|_{{\mathrm{L}}^{2}(\Om)} 
\nonumber .
  \end{align}
By \eqref{2.8} and \eqref{2.11} and by choosing $\e>0$ small, we infer 
\[
\| u \|_{{\mathrm{W}}^{1,p}(\Om)} \, \leq\, C\Big( \| u \|_{{\mathrm{L}}^2(\Om)} +\, \| f \|_{{\mathrm{L}}^{\frac{np}{n+p}}(\Om)} + \| F \|_{{\mathrm{L}}^p(\Om)} +\,\| g \|_{{\mathrm{L}}^p(\p\Om)} \Big).
\]
The estimate \eqref{2.7} follows by combining the above with the ${\mathrm{L}}^{2}$ estimate \eqref{2.5}, together with H\"older inequality and the fact that $\min\big\{p,{np}/(n+p)\big\}>2$. The theorem has been established.
\ep

\smallskip

As a consequence of our result above, we show that the (forward) Robin problem \eqref{1.3} is well posed.

\begin{corollary}[Well-posedness of \eqref{1.3}] \label{th3} Consider \eqref{1.3} and suppose that the coefficients $\A,\B,\K,\L,\M,\xi, s, S$ satisfy \eqref{1.6}-\eqref{3.1}. We further assume that for some $m\geq n$ we have
\beq
\label{3.2}
S \in L ^{\frac{nm}{2n+m}}(\Om;\R^{2}),\ \ \ s \in {\mathrm{L}}^{\frac{m}{2}}(\p\Om;\R^2), \ \ \ \xi \in {\mathrm{L}}^p(\Om), 
\eeq
for some $p>\max\left\{n,{2n}/({n-2}) \right\}$. Then, the problem \eqref{1.3} has a unique weak solution
\[
(u,v) \, \in \, {\mathrm{W}}^{1,\frac{m}{2}}(\Om;\R^{2}) \by {\mathrm{W}}^{1,p}(\Om;\R^{2})
\]
which for any pair of test maps
\[
(\phi,\psi) \, \in \, {\mathrm{W}}^{1,\frac{m}{m-2}} (\Om;\R^{2}) \by {\mathrm{W}}^{1,{\frac{p}{p-1}}} (\Om;\R^2)
\]
it satisfies
\beq \label{3.3}
\ \ \ \left\{
\begin{split}
\!\!\!\! \phantom{\Bigg|} & \int_\Om\Big[ \A :(\D u^\top \D \phi) \, +  \big(\K u-S\big)\cdot \phi \Big] \, \mathrm{d}\mL^n \,+ \int_{\p\Om}\big[ (\ga u-s)\cdot \phi\big] \, \mathrm{d}\mH^{n-1} =0 ,
\\
& \int_\Om\Big[ \B :(\D v^\top \D \psi) \, +  \big(\L v- \xi \M u\big)\cdot \psi \Big] \, \mathrm{d}\mL^n \,+ \int_{\p\Om}\big[ \ga u\cdot \psi\big] \, \mathrm{d}\mH^{n-1} =0 .
\end{split}
\right.
\eeq   
In addition, exists $C>0$ depending only on the data such that
\beq 
\label{3.4}
\left\{
\ \
\begin{split}
   \| u \|_{{\mathrm{W}}^{1,\frac{m}{2}}(\Om)}  & \leq\, C\Big( \| s\|_{{\mathrm{L}}^{\frac{m}{2}}\!(\p\Om)} +\, \| S \|_{L ^{\frac{nm}{2n+m}}\!(\Om)} \Big), \phantom{\Big|}
\\
  \| v \|_{{\mathrm{W}}^{1,p}(\Om)} \, & \leq\, C \| \xi \|_{{\mathrm{L}}^p(\Om)} \, \| u \|_{{\mathrm{L}}^m(\Om)} .
\end{split}
\right.
\eeq
\end{corollary}

\bp By Theorem \ref{th1} applied to the Robin problem \eqref{1.3}(a)-\eqref{1.3}(c) for $p=m/2$ and by recalling that
\[
\frac{m}{2}\, \geq \, \frac{n}{2} \,>\, \frac{2n}{n-2},
\]
for any $s \in {\mathrm{L}}^{\frac{m}{2}}(\p\Om;\R^2)$  and $S \in L ^{\frac{nm}{2n+m}}(\Om;\R^{2})$ exists a unique $u \in {\mathrm{W}}^{1,\frac{m}{2}}(\Om;\R^{2})$ satisfying the $u$-system in \eqref{3.3} for all $\phi \in {\mathrm{W}}^{1,\frac{m}{m-2}} (\Om;\R^{2})$, as well as the estimate \eqref{3.4} for $u$. Fix now $\xi \in {\mathrm{L}}^p(\Om)$. By Theorem \ref{th1} applied to the Robin boundary value problem \eqref{1.3}(b)-\eqref{1.3}(d), it will follow there exists a unique $v \in {\mathrm{W}}^{1,p}(\Om;\R^{2})$ satisfying the $v$-system in \eqref{3.3} for all $\psi \in {\mathrm{W}}^{1,{\frac{p}{p-1}}} (\Om;\R^{2})$, once we have verified that $ \xi \M u \in {\mathrm{L}}^{\frac{np}{n+p}}(\Om)$. By \eqref{2.7} and H\"older's inequality, we estimate
\[
\begin{split}
\| v \|_{{\mathrm{W}}^{1,p}(\Om)} \, & \leq\, C \| \xi \M u\|_{{\mathrm{L}}^{\frac{np}{n+p}}(\Om)}
\\
&\leq\, C \| \xi \|_{{\mathrm{L}}^{\frac{npr}{n+p}}(\Om)}  \| u\|_{{\mathrm{L}}^{\frac{npr'}{n+p}}(\Om)},
\end{split}
\]
for a new $C>0$ and any $r>1$. We select $r:= (n+p)/p$ to obtain
\[
\frac{npr'}{n+p}\, =\, \frac{np}{n+p}\frac{p+n}{p}\, =\, n.
\]
Hence, we conclude that
\[
\begin{split}
\| v \|_{{\mathrm{W}}^{1,p}(\Om)} \, &\leq\, C \| \xi \|_{{\mathrm{L}}^{p}(\Om)}  \| u\|_{{\mathrm{L}}^n(\Om)} 
\, \leq\, C \| \xi \|_{{\mathrm{L}}^{p}(\Om)}  \| u\|_{{\mathrm{L}}^m(\Om)},
\end{split}
\]
because $m\geq n$. The proof is complete. 
\ep

\ms

Now we establish Theorem \ref{th4}. Its proof is a consequence of the next two propositions, utilising the direct method of Calculus of Variations (\cite{D}).

\begin{proposition}[$\I_p$-minimisers] \label{pr4} In the context of Theorem \ref{th4}, the functional $\I_p$ has a constrained minimiser $\big(\vec u_p,\vec v_p ,\xi_p \big) \in \mX^p(\Om)$, as per \eqref{4.15A}.
\end{proposition}

\bp Note that $\mX^p(\Om) \neq \emptyset$, and in fact $\mX^p(\Om)$ is a weakly closed subset of the reflexive Banach space $\mathscr{X}^p(\Om)$ with cardinality greater or equal to that of ${\mathrm{L}}^p(\Om)$. Further, there is an a priori energy bound for the infimum of $\I_p$, in fact uniform in $p$. Indeed, for each $i\in\{1,...,N\}$ let $(u_{0i},v_{0i})$ be the solution to \eqref{1.3} with $\xi \equiv 0$ and $(S_i,s_i)$ as in \eqref{4.2}. Then, by Corollary \ref{th3}, we have $v_{0i} \equiv 0$. Therefore, by \eqref{4.5}-\eqref{4.6} we infer that
\[
\big(\vec u_0,\vec 0 ,0 \big) \, \in \, \mX^p(\Om)
\]
for all $p \in [n,\infty]$, and also, by H\"older inequality and \eqref{4.3}-\eqref{4.7}, we obtain
\[
\begin{split}
\I_p\big(\vec u_0,\vec 0 ,0 \big) \, & \leq\, \frac{N+1}{p} \,+\, E_\infty\big(\vec u_0,\vec 0 ,0 \big)
\, \leq\,  \frac{N+1}{n} \,+\, \sum_{i=1}^N \big\| \tilde v_i\big\|_{{\mathrm{L}}^\infty(\mathrm{B}_i)}
 <\, \infty.
\end{split}
\]
Fix $p$ and consider now a minimising sequence $(\vec u^j,\vec v^{\, j},\xi^{\, j}\,)_{j=1}^\infty$ of $\I_p$ in $\mX^p(\Om)$. Then, for any large enough $j\in\N$ we have 
\[
0\,\leq\, \I_p\big(\vec u^{\,j} ,\vec v^{\, j}, \xi^{\, j}\big) \, \leq\, \frac{N+1}{n} \,+\, \sum_{i=1}^N \big\| \tilde v_i \big\|_{{\mathrm{L}}^\infty(\mathrm{B}_i)}\, +\,1.
\]
By Corollary \ref{th3}, we have the estimates
\beq
\label{4.13}
\left\{
\ \
\begin{split}
\phantom{_\big|}  \big\| \vec u^{\,j} \big\|_{{\mathrm{W}}^{1,\frac{m}{2}}(\Om)}   \, & \leq\, C \, \max_{i=1,...,N} \Big(\| S_i \|_{L ^{\frac{nm}{2n+m}}\!(\Om)} +\,\| s_i\|_{{\mathrm{L}}^{\frac{m}{2}}\!(\p\Om)}\Big), \phantom{\Big|}
\\
\big\| \vec v^{\, j}\big\|_{{\mathrm{W}}^{1,p}(\Om)} \, & \leq\, C \, \| \xi^{\, j}\|_{{\mathrm{L}}^p(\Om)} \, \big\| \vec u^{\, j}\big\|_{{\mathrm{W}}^{1,\frac{m}{2}}(\Om)} .
\end{split}
\right.
\eeq
If \eqref{4.11} is satisfied, then by all the above and \eqref{4.4} we have the estimate
\[
\| \xi^{\,j} \|_{{\mathrm{L}}^p(\Om)} \, \leq\, \| \xi^{\,j} \|_{\dot{\mathrm{L}}^p(\Om)} \, \leq\, \frac{1}{\al} \bigg( \frac{N+1}{n}+ \sum_{i=1}^N \big\| \tilde v_i\big\|_{{\mathrm{L}}^\infty(\mathrm{B}_i)} +\, 1\bigg).
\]
If on the other hand \eqref{4.12} is satisfied, then we readily have
\[
\| \xi^{\, j}\|_{{\mathrm{L}}^p(\Om)} \, \leq\, \| \xi^{\, j} \|_{{\mathrm{L}}^\infty(\Om)} \, \leq\, M .
\]
Hence, in both cases of either \eqref{4.11} or \eqref{4.12}, we have the uniform bound
\beq
\label{4.14}
\sup_{j\in\N}\, \| \xi^{\, j}\|_{{\mathrm{L}}^p(\Om)} \, <\, \infty .
\eeq
By the estimates \eqref{4.13}-\eqref{4.14} and standard weak and strong compactness arguments, there exists a weak limit in the Banach space
\[
(\vec u_p,\vec v_p , \xi_p \big) \, \in \mathscr{X}^p(\Om)
\]
and a subsequence $(j_k)_1^\infty$ such that as $j_k \to \infty$ we have
\[
\left\{ \ \
\begin{array}{ll}
\xi^{\, j} \weak \xi_p, & \text{ in }{\mathrm{L}}^p(\Om),
\ms
\\
\vec u^{\, j} \weak \vec u_p, & \text{ in }{\mathrm{W}}^{1,\frac{m}{2}}(\Om;\R^{2\by N}), \ms
\\
\vec u^{\, j} \larrow \vec u_p, & \text{ in }{\mathrm{L}}^{\frac{m}{2}}(\Om;\R^{2\by N}), \ms
\\
\vec v^{\, j} \weak \vec v_p, & \text{ in }{\mathrm{W}}^{1,p}(\Om;\R^{2\by N}), \ms
\\
\vec v^{\, j} \larrow \vec v_p, & \text{ in } {\mathrm{C}}^0(\overline{\Om};\R^{2\by N}).
\end{array}
\right.
\]
Note that in this paper we utilise the standard practice of passing to subsequences as needed perhaps without explicit mention. To show that $(\vec u_p,\vec v_p , \xi_p \big)$ lies in $\mX^p(\Om)$, we argue as follows. First we show that for any $M \in[0,\infty]$ the constraint
\[
0 \, \leq\, \xi^{\,j}   \leq\, M \ \text{ a.e. on }\Om,
\]
is weakly closed in ${\mathrm{L}}^p(\Om)$, namely
\beq
\label{4.15}
{\mathrm{L}}^p(\Om;[0,M])\,=\, \Big\{\eta \in {\mathrm{L}}^p(\Om):\, \ 0 \, \leq\, \eta\,   \leq\, M \ \text{ a.e. on }\Om\Big\}
\eeq
is weakly closed. To this aim, let $\xi^{\, j} \weak \xi_p$ in ${\mathrm{L}}^p(\Om)$ as $j_k\to \infty$. Then, for any measurable set $E\sub \Om$ with positive measure $\mL^n(E)>0$, by integrating the last inequality over $E$, the averages satisfy
\[
0 \, \leq\, \av_E \xi^{\,j}\, \mathrm d \mL^n   \leq\, M 
\]
and therefore
\[
\ \ \av_E \xi_p\, \mathrm d \mL^n =\lim_{j_k\to \infty} \left( \, \av_E \xi^{\,j}\, \mathrm d \mL^n \! \right)  \in  [0,M]. 
\]
By selecting $E:= \mB_\rho(x)$ for $x\in\Om$ and $\rho \in (0,\dist(x,\p\Om))$, the Lebesgue differentiation theorem allows us to infer
\[
\ \ \xi_p(x) \, =\lim_{\rho\to 0} \left(\, \av_{\mB_\rho(x)} \xi_p\, \mathrm d \mL^n \! \right)   \in  [0,M], \ \ \text{ for a.e. }x\in \Om. 
\]
To conclude that $(\vec u_p,\vec v_p , \xi_p \big) \in \mX^p(\Om)$, we must pass to the weak limit in the equations $(a)_i-(d)_i$ in \eqref{4.8}. The only convergence that needs to be justified that of the nonlinear source term $\xi \M u_i$ in $(b)_i$. To this end, note that by our assumption $p>\frac{2n}{n-2}$, we have the inequality
\[
{\frac{p}{p-1}}\, <\, \frac{n}{2} \, \leq \, \frac{m}{2}. 
\]
Thus, since $u_i^{\, j} \larrow u_{pi}$ in ${\mathrm{L}}^{{\frac{m}{2}}}(\Om;\R^{2})$ as $j_k \to \infty$, we have that
\[
 u_i^{\, j} \larrow u_{pi} \ \text{ in } {\mathrm{L}}^{{\frac{p}{p-1}}}(\Om;\R^{2})
\]
as $j_k \to \infty$. Hence, since $\xi^{\, j} \weak \xi_{p}$ in ${\mathrm{L}}^{p}(\Om)$, it follows that 
\[
\int_\Om \big(\xi^{\, j}\, \M u^{\, j}_i \big) \cdot \phi \, \mathrm d \mL^n \larrow \int_\Om \big(\xi_p \, \M u_{pi} \big) \cdot \phi \,  \mathrm d \mL^n
\]
for any $\phi \in {\mathrm{C}}^\infty_c(\Om;\R^2)$ as $j_k \to \infty$, as a consequence of the weak-strong continuity of the duality pairing between ${\mathrm{L}}^p(\Om)$ and ${\mathrm{L}}^{{\frac{p}{p-1}}}(\Om)$. In conclusion, we infer that indeed $(\vec u_p,\vec v_p , \xi_p \big) \in \mX^p(\Om)$ by passing to the weak limit in the equations $(a)_i-(d)_i$ defining \eqref{4.8}.

We now show that $(\vec u_p,\vec v_p , \xi_p \big) \in \mX^p(\Om)$ is indeed a minimiser of $\I_p$. To this aim, note that for any $\al \in [0,\infty)$ the functional $\al \| \cdot \|_{\dot{\mathrm{L}}^p(\Om)}$ is convex and strongly continuous on the reflexive space ${\mathrm{L}}^p(\Om)$, by \eqref{4.6}-\eqref{4.7}. Thus, it is weakly lower semi-continuous. Similarly, note that for each index $i\in\{1,...,N\}$ the functional $\|\cdot - \, \tilde v_i\|_{\dot{\mathrm{L}}^p(\mathrm{B}_i)}$ is strongly continuous on ${\mathrm{L}}^p(\mathrm{B}_i)$. In conclusion, we have
\[
\begin{split}
\I_p\big(\vec u_p,\vec v_p ,\xi_p \big) \,&=\, \al \| \xi_p \|_{\dot{\mathrm{L}}^p(\Om)}  +\, \sum_{i=1}^N \big\|v_{pi} - \tilde v_i\big\|_{\dot{\mathrm{L}}^p(\mathrm{B}_i)} 
\\
&\leq \liminf_{j_k\to \infty} \left\{ \al \| \xi^{\, j} \|_{\dot{\mathrm{L}}^p(\Om)}  +\, \sum_{i=1}^N \big\|v^{\, j}_{i} - \tilde v_i\big\|_{\dot{\mathrm{L}}^p(\mathrm{B}_i)}  \right\}
\\
&= \liminf_{j_k\to \infty} \, \I_p\big(\vec u^{\,j} ,\vec v^{\, j}, \xi^{\, j}\big)
\\
&=\, \inf\Big\{\I_p(\vec u,\vec v,\xi \big) :\ (\vec u,\vec v ,\xi \big) \in \mX^p(\Om) \Big\}.
\end{split}
\]
The proposition ensues.
\ep

Our next result below establishes the existence of minimisers for ${\mathrm{I}}_\infty$ and the approximation by minimisers of the ${\mathrm{I}}^p$ functionals as $p\to\infty$, therefore completing the proof of Theorem \ref{th4}.

\begin{proposition}[$\I_\infty$-minimisers] \label{pr5} The functional $\I_\infty$ given by \eqref{4.5} has a constrained minimiser $(\vec u_\infty,\vec v_\infty ,\xi_\infty) \in \mX^\infty(\Om)$, as in \eqref{4.15B}. Additionally, exists a subsequence $(p_j)_1^\infty$ such that the $\I_{p_j}$-minimisers $\big(\vec u_{p_j},\vec v_{p_j} ,\xi_{p_j} \big)$ satisfy \eqref{4.16} as $j\to\infty$. 
\end{proposition}

\bp We continue from the proof of Proposition \ref{pr4}. The map $(\vec u_0,\vec 0,0)$ constructed therein provides an energy bound uniform in $p$ and in view of \eqref{4.8}-\eqref{4.9} we also have $(\vec u_0,\vec 0,0) \in\mX^\infty(\Om)$. Fix $q>n$ and $p\geq q$. By H\"older inequality and minimality, we estimate
\[
\begin{split}
\I_q\big(\vec u_{p},\vec v_{p} ,\xi_{p} \big)\, &\leq\, \I_p\big(\vec u_{p},\vec v_{p} ,\xi_{p} \big)
\\
 & \leq\, \I_p(\vec u_0,\vec 0,0)
 \\
 & \leq\, \I_\infty(\vec u_0,\vec 0,0) \,+\, \frac{N+1}{n}
  \\
  & \leq \, \sum_{i=1}^N \big\| \tilde v_i\big\|_{{\mathrm{L}}^\infty(\mathrm{B}_i)} \,+\, \frac{N+1}{n},
\end{split}
\]
which is uniform in $p$. If \eqref{4.11} is satisfied, then by the above estimate we have
\[
\| \xi_p \|_{{\mathrm{L}}^q(\Om)} \, \leq\, \| \xi_p \|_{\dot{\mathrm{L}}^p(\Om)} \, \leq\, \frac{1}{\al} \bigg(\sum_{i=1}^N \big\| \tilde v_i\big\|_{{\mathrm{L}}^\infty(\mathrm{B}_i)} +\, \frac{N+1}{n} \bigg).
\]
If on the other hand \eqref{4.12} is satisfied, then we immediately have
\[
\| \xi_p\|_{{\mathrm{L}}^q(\Om)}  \leq\, \| \xi_p \|_{{\mathrm{L}}^\infty(\Om)} \leq\, M .
\]
Hence, in both cases of either \eqref{4.11} or \eqref{4.12}, we have the uniform bound
\beq
\label{4.14}
\sup_{p\geq q}\, \| \xi_p\|_{{\mathrm{L}}^q(\Om)} \, <\, \infty .
\eeq
By the above estimates, Corollary \ref{th3} (see \eqref{4.13}) and standard compactness arguments yield that there exists a subsequence $(p_j)_1^\infty$ and a limit
\[
(\vec u_\infty,\vec v_\infty , \xi_\infty \big) \, \in \bigcap_{n<q<\infty}\mW^q(\Om)
\]
such that \eqref{4.16} holds true as $j \to \infty$. In addition, under either assumptions \eqref{4.11} or \eqref{4.12} we have that $\xi_\infty \in {\mathrm{L}}^\infty(\Om)$. Indeed, under \eqref{4.11} by integrating the constraint $0\leq \xi_p \leq M$ on a measurable set $E\sub \Om$ of positive measure we have
\[
0 \, \leq\, \av_E \xi_p\, \mathrm d \mL^n   \leq\, M 
\]
and therefore
\[
\ \ \av_E \xi_\infty\, \mathrm d \mL^n =\lim_{p_j\to \infty} \left( \, \av_E \xi_p\, \mathrm d \mL^n \! \right)  \in  [0,M]. 
\]
Then, for $E:= \mB_\rho(x)$, $x\in\Om$ and $\rho \in (0,\dist(x,\p\Om))$ we deduce
\[
\ \ \xi_\infty(x) \, =\lim_{\rho\to 0} \left( \, \av_{\mB_\rho(x)} \xi_\infty\, \mathrm d \mL^n \! \right)   \in  [0,M], \ \ \text{ for a.e. }x\in \Om, 
\]
for any $M>0$. If on the other hand \eqref{4.11} is satisfied, then by the weak lower-semicontinuity of the functional $\|\cdot\|_{\dot{\mathrm{L}}^q(\Om)}$ on ${\mathrm{L}}^q(\Om)$, we have
\[
\begin{split}
\|\xi_\infty \|_{{\mathrm{L}}^\infty (\Om)} \,&=\, \lim_{q\to \infty}\| \xi_\infty \|_{\dot{\mathrm{L}}^q(\Om)} 
\\
&\leq\, \liminf_{q\to \infty} \left(\liminf_{p_j\to \infty} \| \xi_p \|_{\dot{\mathrm{L}}^q(\Om)}\right)
\\
&\leq\, \liminf_{q\to \infty} \left(\liminf_{p_j\to \infty} \frac{1}{\al} \sum_{i=1}^N \big\| \tilde v_i\big\|_{{\mathrm{L}}^\infty(\mathrm{B}_i)} \right) 
\\
&= \frac{1}{\al} \sum_{i=1}^N \big\| \tilde v_i\big\|_{{\mathrm{L}}^\infty(\mathrm{B}_i)} .
\end{split}
\]
Further, by passing to the limit as $p_j\to \infty$ in $(a)_i-(d)_i$ of \eqref{4.8} as in the proof of Proposition \ref{pr4}, we see that the limit $(\vec u_\infty,\vec v_\infty , \xi_\infty \big)$ lies in $\mX^\infty(\Om)$. It remains to prove that $(\vec u_\infty,\vec v_\infty , \xi_\infty \big)$ is a minimiser of $\I_\infty$ and that the energies converge. Fix an arbitrary $(\vec u,\vec v, \xi \big) \in \mX^\infty(\Om)$. Since $p_j\geq q$ for large $j\in\N$, by minimality we have
\[
\begin{split}
\I_\infty\big(\vec u_\infty,\vec v_\infty ,\xi_\infty \big) \,&=\, \lim_{q\to \infty} \I_q\big(\vec u_\infty,\vec v_\infty ,\xi_\infty \big)
\\
&\leq\, \liminf_{q\to \infty} \Big(\liminf_{p_j\to \infty}\, \I_q\big(\vec u_p,\vec v_p ,\xi_p \big)  \Big)
\\
&\leq\, \liminf_{p_j\to \infty} \I_p\big(\vec u_p,\vec v_p ,\xi_p \big) 
\\
&\leq\, \limsup_{p_j\to \infty} \I_p\big(\vec u_p,\vec v_p ,\xi_p \big)  
\\
&\leq\, \limsup_{p_j\to \infty} \I_p\big(\vec u,\vec v,\xi \big)  
\\
&=\,  \I_\infty\big(\vec u,\vec v ,\xi \big) , 
\end{split}
\]
for any $(\vec u,\vec v, \xi \big) \in \mX^\infty(\Om)$. Therefore, $\big(\vec u_\infty,\vec v_\infty ,\xi_\infty \big)$ is a minimiser of $\I_\infty$ in $\mX^\infty(\Om)$. The choice $(\vec u,\vec v, \xi \big):=\big(\vec u_\infty,\vec v_\infty ,\xi_\infty \big)$ in the above inequality implies
\[
\I_\infty\big(\vec u_\infty,\vec v_\infty ,\xi_\infty \big) \larrow \, \I_p\big(\vec u_p,\vec v_p ,\xi_p \big),
\]
as $j\to\infty$. The proof of Proposition \ref{pr5} is complete.
\ep

Before proving Theorem \ref{th7}, we use it to obtain an additional piece of information on the variational inequality \eqref{5.20}.

\begin{corollary} In the setting of Theorem \ref{th7}, in the case under assumption \eqref{4.11} (where $M=\infty$ and $\al>0$) the variational inequality \eqref{5.20} for the constrained minimiser implies the next test-function-free relations:
\begin{align}
\sum_{i=1}^N  \big(\M  u_{pi}  \big)  \cdot \psi_{pi}  \, & \geq\, 0, \ \ \text{ a.e. on }\{\xi_p=0\},
\label{5.23}
\\
\al\,\frac{ |\xi_p|_{(p)}^{p-2} \, \xi_p }{ \mL^{n}(\Om)  \|\xi_p \|_{\dot{\mathrm{L}}^p(\Om)}^{p-1}}   \, + \,  \sum_{i=1}^N  \big(\M  u_{pi}  \big)  \cdot \psi_{pi}    \,& =\, 0,  \ \ \text{ a.e. on }\{\xi_p>0\}.
\label{5.24}
\end{align}
\end{corollary}

\bp To see \eqref{5.23}, note that if $M=\infty$, then by testing in \eqref{5.20} against $\eta := \xi_p +\theta$ where $\theta \in {\mathrm{L}}^p(\Om;[0,\infty))$, we obtain 
\[
\int_\Om \theta \left( \al\, \frac{\mathrm d [\mu_p(\xi_p)]}{\mathrm{d}\mL^n} \, + \sum_{i=1}^N  \big(\M  u_{pi}  \big)  \cdot \psi_{pi} \right) \, \mathrm{d}\mL^n \, \geq\, 0,
\]
for any $\theta \in {\mathrm{L}}^p(\Om,[0,+\infty))$, which yields
\beq
\label{5.25}
\al\,\frac{ |\xi_p|_{(p)}^{p-2} \, \xi_p }{ \mL^{n}(\Om)  \|\xi_p \|_{\dot{\mathrm{L}}^p(\Om)}^{p-1}}   \, + \,  \sum_{i=1}^N  \big(\M  u_{pi}  \big)  \cdot \psi_{pi}    \, \geq \, 0,  \ \ \text{ a.e. on }\Om.
\eeq
From the above inequality we readily deduce \eqref{5.23}. To see \eqref{5.24}, we fix a point $x\in \{\xi_p>0\}$, $t>0$ small and $\rho \in(0,\dist(x,\p\Om)$ and test against the function
\[
\eta \,:=\, \xi_p -t \chi_{\{\xi_p>t\}\cap \mB_\rho(x)} \, \in {\mathrm{L}}^p(\Om;[0,\infty).
\]
Then, by \eqref{5.20} we get
\[
t \int_{\mB_\rho(x)} \chi_{\{\xi_p>t\}} \left( \al\, \frac{\mathrm d [\mu_p(\xi_p)]}{\mathrm{d}\mL^n} \, + \sum_{i=1}^N  \big(\M  u_{pi}  \big)  \cdot \psi_{pi} \right) \, \mathrm{d}\mL^n \, \leq\, 0,
\]
which by diving by $t\mL^n(\mB_\rho(x))$, letting $t\to0$, using the Dominated Convergence theorem and letting $\rho\to0$ yields 
\[
\lim_{\rho\to0} \ \av_{\mB_\rho(x)} \chi_{\{\xi_p>0\}} \left( \al\, \frac{\mathrm d [\mu_p(\xi_p)]}{\mathrm{d}\mL^n} \, + \sum_{i=1}^N  \big(\M  u_{pi}  \big)  \cdot \psi_{pi} \right) \, \mathrm{d}\mL^n \, \leq\, 0.
\]
Now, \eqref{5.24} follows as a consequence of the Lebesgue differentiation theorem and \eqref{5.25}. The proof is complete.
\ep

The proof of Theorem \ref{th7} consists of a few sub-results. We begin by computing the derivative of $\I_p$.

\begin{lemma} \label{pr7} The functional $\I_p : \mathscr{X}^p(\Om) \larrow \R$ is Frech\'et differentiable and its derivative
\[
\mathrm d \hspace{0.5pt} \I_p \ : \ \ \mathscr{X}^p(\Om) \larrow \big(\mathscr{X}^p(\Om)\big)^*
\]
which maps 
\[
(\vec u, \vec v,\xi) \mapsto \big(\mathrm d \hspace{0.5pt} \I_p\big)_{ (\vec u, \vec v,\xi)}
\]
is given for all $(\vec u, \vec v,\xi), (\vec z, \vec w,\eta) \in \mathscr{X}^p(\Om)$ by the formula
\beq 
\label{5.1}
\big(\mathrm d \hspace{0.5pt} \I_p\big)_{ (\vec u, \vec v,\xi)} (\vec z, \vec w,\eta)  \,=\, p\int_{\p\Om} \vec w : \mathrm d [\vec \nu_p(\vec v)]\, + \, \al p \int_\Om  \eta\, \mathrm d [\mu_p(\xi)].
\eeq 
\end{lemma}

\bp The Frech\'et differentiability of $\I_p$ follows from well-known results on the differentiability of norms on Banach spaces and our $p$-regularisations in \eqref{4.6}-\eqref{4.7}. To compute the Frech\'et derivative, we use directional differentiation. For any fixed $(\vec u, \vec v,\xi), (\vec z, \vec w,\eta) \in \mathscr{X}^p(\Om)$, we have
\[
\begin{split}
\frac{1}{p}\big(\mathrm d \hspace{0.5pt} \I_p\big)_{ (\vec u, \vec v,\xi)} (\vec z, \vec w,\eta)  \, & =\, \frac{1}{p} \frac{\mathrm d}{\mathrm d \e}\Big|_{\e=0}\I_p \Big((\vec u, \vec v,\xi)+\e (\vec z, \vec w,\eta)\Big)
\\
&=\,\frac{1}{p} \sum_{i=1}^N \frac{\mathrm d}{\mathrm d \e}\Big|_{\e=0}\left( \, \av_{\mathrm{B}_i}\big(\big| v_i+\e w_i -\tilde v_i \big|_{(p)} \big)^p\, \mathrm d \mH^{n-1} \! \right)^{\! \frac{1}{p}}
\\
& \ \ \ \ + \, \al \, \frac{1}{p} \frac{\mathrm d}{\mathrm d \e}\Big|_{\e=0}\left( \, \av_{\Om}\big(|\xi+\e \eta|_{(p)}\big)^p\, \mathrm d \mL^{n} \! \right)^{\! \frac{1}{p}}
\\
& = \, \sum_{i=1}^N  \left( \, \av_{\mathrm{B}_i}\big(\big| v_i -\tilde v_i \big|_{(p)} \big)^p\, \mathrm d \mH^{n-1} \! \right)^{\! \frac{1}{p}-1} \centerdot
\\
&\ \ \ \ \,  \centerdot \av_{\mathrm{B}_i}\big(|v_i-\tilde v_i|_{(p)}\big)^{p-2} (v_i-\tilde v_i) \cdot w_i\, \mathrm d \mH^{n-1} 
\\
& \ \ \ \ + \al \left( \, \av_{\Om}\big(|\xi|_{(p)}\big)^p\, \mathrm d \mL^{n} \! \right)^{\! \frac{1}{p} -1}
\av_{\Om}(|\xi|_{(p)})^{p-2}\xi \,\eta \, \mathrm d \mL^{n}.
\end{split}
\]
Hence, \eqref{5.1} follows in view of \eqref{5.3}-\eqref{5.2}. The lemma ensues.
\ep

In order to derive the variational inequality that any $p$-minimiser as in \eqref{4.15A} satisfies, we define a mapping which expresses the PDE constraints of the admissible class in \eqref{4.8}. Thus, we set
\beq \label{5.6}
\J \ : \ \ \mathscr{X}^p(\Om) \larrow \Big[\big( {\mathrm{W}}^{1,\frac{m}{m-2}} (\Om;\R^{2})\big)^*  \by \big({\mathrm{W}}^{1,\frac{p}{p-1}} (\Om;\R^{2})\big)^*\Big]^N,
\eeq
\beq  \label{5.7}
\begin{split}
& \left\langle \J(\vec u, \vec v, \xi )  ,\, (\vec \phi, \vec \psi\, )\right\rangle\, := \, \left[
\begin{array}{c}
 \Big\langle \J^1_1(\vec u, \vec v, \xi )  ,\,  \phi_1 \Big\rangle \ms
 \\
  \Big\langle \J^2_1(\vec u, \vec v, \xi )  ,\,  \psi_1 \Big\rangle 
 \\
 \vdots
 \\
 \Big\langle \J^1_N(\vec u, \vec v, \xi )  ,\, \phi_N \Big\rangle \ms
 \\
 \Big\langle \J^2_N(\vec u, \vec v, \xi )  ,\,  \psi_N\Big\rangle 
 \end{array}
\right] \, \in \, \R^{2N},
\end{split}
\eeq   
where, for any $j\in\{1,2\}$, and $i\in\{1,...,N\}$ and any test maps 
\[
(\phi_i, \psi_i ) \ \in {\mathrm{W}}^{1,\frac{m}{m-2}} (\Om;\R^{2}) \by {\mathrm{W}}^{1,\frac{p}{p-1}} (\Om;\R^{2}), 
\]
the component $\J^j_i$ of $\J$ is given by
\beq 
 \label{5.8}
\left\{\  
\begin{split}
 \Big\langle \J^1_i(\vec u, \vec v, \xi )  ,\, \phi_i \Big\rangle\, := & \int_\Om\Big[ \A :(\D u_i^\top \D \phi_i) \, +  \big(\K u_i-S_i\big)\cdot \phi_i \Big] \, \mathrm{d}\mL^n
\\
&  + \int_{\p\Om}\big[ (\ga u_i-s_i)\cdot \phi_i\big] \, \mathrm{d}\mH^{n-1}
\end{split}
\right.
\eeq 
\beq 
 \label{5.9}
\left\{\  
\begin{split}
 \Big\langle \J^2_i(\vec u, \vec v, \xi )  ,\, \psi_i \Big\rangle\, := & \int_\Om\Big[ \B :(\D v_i^\top \D \psi_i) \, +  \big(\L v_i- \xi \M u_i\big)\cdot \psi_i \Big] \, \mathrm{d}\mL^n 
\\
&   + \int_{\p\Om}\big[ \ga v_i\cdot \psi_i\big] \, \mathrm{d}\mH^{n-1}.
\end{split}
\right.
\eeq   
Let us also define for any $M\in[0,\infty]$ the following weakly closed convex subset of the Banach space $\mathscr{X}^p(\Om)$:
\beq
 \label{5.10}
\mathscr{X}^p_M(\Om) \,:=\, {\mathrm{W}}^{1,\frac{m}{2}}(\Om;\R^{2\by N}) \by {\mathrm{W}}^{1,p}(\Om;\R^{2\by N}) \by {\mathrm{L}}^p(\Om;[0,M]).
\eeq 
Then, in view of \eqref{5.6}-\eqref{5.10}, we may reformulate the admissible class $\mX^p(\Om)$ of the minimisation problem \eqref{4.15A} as
\beq
\label{5.11}
\mX^p(\Om) \, =\, \Big\{ \big(\vec u, \vec v,\xi\big) \in \mathscr{X}^p_M(\Om)\, : \ \J \big(\vec u, \vec v,\xi\big)=0\Big\}.
\eeq

We now compute the derivative of $\J$ above and prove that it is a $C^1$ submersion.

\begin{lemma} \label{le8}
The map $\J$ defined by \eqref{5.6}-\eqref{5.10} is a continuously differentiable submersion and its Frech\'et derivative
\beq \label{5.12}
\begin{split}
\mathrm d \hspace{0.5pt} \J  \ : \ \ \mathscr{X}^p(\Om) \larrow \mathscr{L} & \left(\mathscr{X}^p(\Om), \Big[\big( {\mathrm{W}}^{1,\frac{m}{m-2}} (\Om;\R^{2})\big)^*  \by \big({\mathrm{W}}^{1,\frac{p}{p-1}} (\Om;\R^{2})\big)^*\Big]^N \right),
 \end{split}
\eeq
which maps
\[
(\vec u, \vec v,\xi) \mapsto \big(\mathrm d \hspace{0.5pt} \J\big)_{ (\vec u, \vec v,\xi)}
\]
is given by
\beq  \label{5.13}
\begin{split}
& \left\langle \big(\mathrm d \hspace{0.5pt} \J\big)_{ (\vec u, \vec v,\xi)} (\vec z, \vec w,\eta) ,\, (\vec \phi, \vec \psi\, )\right\rangle\, = \, \left[
\begin{array}{c}
 \Big\langle \big(\mathrm d \hspace{0.5pt} \J^1_1\big)_{ (\vec u, \vec v,\xi)} (\vec z, \vec w,\eta)  ,\,  \phi_1 \Big\rangle \ms
 \\
  \Big\langle \big(\mathrm d \hspace{0.5pt} \J^2_1\big)_{ (\vec u, \vec v,\xi)} (\vec z, \vec w,\eta)  ,\,  \psi_1  \Big\rangle 
 \\
 \vdots
 \\
 \Big\langle \big(\mathrm d \hspace{0.5pt} \J^1_N\big)_{ (\vec u, \vec v,\xi)} (\vec z, \vec w,\eta)   ,\,  \phi_N \Big\rangle \ms
 \\
 \Big\langle \big(\mathrm d \hspace{0.5pt} \J^2_N\big)_{ (\vec u, \vec v,\xi)} (\vec z, \vec w,\eta)  ,\, \psi_N \Big\rangle 
 \end{array}
\right] .
\end{split}
\eeq   
In \eqref{5.13}, for each $i \in\{1,...,N\}$ and $j\in\{1,2\}$, the component $\big(\mathrm d \hspace{0.5pt} \J^1_N\big)_{ (\vec u, \vec v,\xi)}$ of the derivative is given for any test functions 
\[
(\vec \phi, \vec \psi\, ) \ \in {\mathrm{W}}^{1,\frac{m}{m-2}} (\Om;\R^{2\by N}) \by {\mathrm{W}}^{1,\frac{p}{p-1}} (\Om;\R^{2\by N})
\]
by the expressions
\beq 
 \label{5.14}
\left\{\!
\begin{split}
 \left\langle \big(\mathrm d \hspace{0.5pt} {\J^1_i}\big)_{ (\vec u, \vec v,\xi)} (\vec z, \vec w,\eta) ,\,  \phi_i \right\rangle\, = & \int_\Om\Big[ \A :(\D z_i^\top \D \phi_i) \, + (\K z_i)\cdot \phi_i \Big] \, \mathrm{d}\mL^n
\\
& \ + \int_{\p\Om}(\ga z_i)\cdot \phi_i \, \mathrm{d}\mH^{n-1},
\end{split}
\right.
\eeq   
\beq 
 \label{5.15}
\left\{ \!
\begin{split}
 \Big\langle \big(\mathrm d \hspace{0.5pt} {\J^2_i} & \big)_{ (\vec u, \vec v,\xi)}  (\vec z, \vec w,\eta)  ,\,  \psi_i  \Big\rangle  \, =  \int_\Om\Big[ \B :(\D w_i^\top \D \psi_i) 
 \\
 & +  \Big(\L w_i -\M (\eta u_i +\xi z_i )\Big)\cdot \psi_i \Big] \, \mathrm{d}\mL^n
\, + \int_{\p\Om} (\ga w_i)\cdot \psi_i \, \mathrm{d}\mH^{n-1},
\end{split}
\right.
\eeq   
and any $(\vec u, \vec v, \xi ), (\vec z, \vec w, \eta  ) \in \mathscr{X}^p(\Om)$. 
\end{lemma}

\bp The mapping $\J$ is bounded linear in all arguments except for the products $\xi \M u_i$ appearing in the component $\J^2_i$. Therefore, it is quadratic and continuously Frech\'et differentiable in the space $\mathscr{X}^p(\Om)$. The Frech\'et derivative can be computed by using directional differentiation
\[
 \big(\mathrm d \hspace{0.5pt} \J\big)_{ (\vec u, \vec v,\xi)}(\vec z, \vec w, \eta ) \, = \, \frac{\mathrm d}{\mathrm d \e}\Big|_{\e=0} \J\Big((\vec u, \vec v, \xi )+\e (\vec z, \vec w, \eta )\Big) 
\]
and the observation that all the terms are affine except for the quadratic term $\M (\xi+\e\eta) (u_i+\e z_i)$ whose derivative as $\e=0$ is $\M (\eta u_i +\xi z_i )$. To conclude, we must show that $\J$ is a submersion, namely for any $(\vec u, \vec v, \xi ) \in \mathscr{X}^p(\Om)$, the differential at this point 
\[
\big(\mathrm d \hspace{0.5pt} \J\big)_{ (\vec u, \vec v,\xi)} \ : \ \ \mathscr{X}^p(\Om) \larrow \Big[\big( {\mathrm{W}}^{1,\frac{m}{m-2}} (\Om;\R^{2})\big)^*  \by \big({\mathrm{W}}^{1,\frac{p}{p-1}} (\Om;\R^{2})\big)^*\Big]^N
\]
is a surjective linear map. To this aim, for each $i\in\{1,...,N\}$ fix a pair of functionals
\[
(\Phi_i, \Psi_i ) \ \in \big( {\mathrm{W}}^{1,\frac{m}{m-2}} (\Om;\R^{2})\big)^*  \by \big({\mathrm{W}}^{1,\frac{p}{p-1}} (\Om;\R^{2})\big)^*.
\]
By well-known results (see e.g.\ \cite{Ad}), it follows that exist
\[
\left\{ \ \
\begin{split}
& (f_i,F_i) \, \in {\mathrm{L}}^{\frac{m}{2}} (\Om;\R^{2})  \by {\mathrm{L}}^{\frac{m}{2}} (\Om;\R^{2\by n}),
\\
&(g_i,G_i) \, \in {\mathrm{L}}^{p} (\Om;\R^{2})  \by {\mathrm{L}}^{p} (\Om;\R^{2\by n}),
\end{split}
\right.
\]
such that for any $\phi_i \in {\mathrm{W}}^{1,\frac{m}{m-2}} (\Om;\R^{2})$ and $\psi_i \in {\mathrm{W}}^{1,\frac{p}{p-1}} (\Om;\R^{2})$, the next representation formulas hold true
\beq
\label{5.16}
\left\{\ \ 
\begin{split}
\langle \Phi_i, \phi_i \rangle \, &=\, \int_\Om \big(f_i\cdot \phi_i \,+\, F_i :\D \phi_i\big)\, \mathrm d\mL^n, 
\\
 \langle \Psi_i, \psi_i \rangle \, &=\, \int_\Om \big(g_i\cdot \psi_i \,+\, G_i :\D \psi_i\big)\, \mathrm d\mL^n.
 \end{split}
 \right.
\eeq
By \eqref{5.12}-\eqref{5.16}, the surjectivity of $\J'(\vec u, \vec v, \xi )$ is equivalent to the weak solvability in $(z_i,w_i)$ of the systems
\[
\left\{\ \ 
\begin{array}{ll}
  -\div(\D z_i \hspace{1pt} \A )\,+\, \K z_i\, =\, f_i-\div F_i, & \ \ \text{ in }\Om, \ms
\\
\ \  (\D z_i \hspace{1pt} \A - F_i)\hspace{1pt} \mathrm n\,+\, \ga z_i\, =\, 0, & \ \ \text{ on }\p\Om,
\end{array}
\right.
\]
and
\[
\left\{\ \ 
\begin{array}{ll}
  -\div( \D w_i \hspace{1pt} \B)\,+\, \L w_i\, =\, \big( g_i+\M(\eta u_i +\xi z_i) \big)-\div \,G_i, & \ \ \text{ in }\Om, \ms
\\
\ \ (\D w_i \hspace{1pt} \B - G_i)\hspace{1pt}\mathrm n\,+\, \ga w_i\, =\, 0, & \ \ \text{ on }\p\Om,
\end{array}
\right.
\]
for all indices $i\in\{1,...,N\}$, some $\eta \in {\mathrm{L}}^p(\Om)$ and with $\A,\B,\K,\L,\M,u_i,\xi,\ga,f_i,F_i$, $g_i,G_i$ being fixed coefficients. The solvability of the above systems follows from Theorem \ref{th1} and Corollary \ref{th3}. The proof is complete.
\ep

Now we derive the variational inequality in $\mathrm L^p$ by employing the Kuhn-Tucker theory of generalised Lagrange multipliers.

\begin{proposition}[The variational inequality] \label{pr9}
For any $p>2n/(n-2)$, there exist Lagrange multipliers
\[
\big(\vec \phi_p, \vec \psi_p \big) \ \in {\mathrm{W}}^{1,\frac{m}{m-2}} (\Om;\R^{2\by N}) \by {\mathrm{W}}^{1,\frac{p}{p-1}} (\Om;\R^{2\by N})
\]
associated with the minimisation problem \eqref{4.15A}, such that the constrained minimiser $\big(\vec u_p, \vec v_p, \xi_p\big) \in \mX^p(\Om)$ satisfies for any $(\vec z, \vec w,\eta)$ in the convex set $\mathscr{X}^p_M(\Om)$ that
\beq
\label{5.17}
\begin{split}
  \frac{1}{p} \big(\mathrm d \hspace{0.5pt} \I_p\big)_{ (\vec u_p, \vec v_p,\xi_p)}   \big(\vec z, \vec w, \eta-\xi_p\big) 
 \, & \geq  \, \sum_{i=1}^N \bigg\langle \big(\mathrm d \hspace{0.5pt} \J^1_i\big)_{ (\vec u_p, \vec v_p,\xi_p)}\big(\vec z, \vec w,\eta-\xi_p\big),\, \phi_{pi}   \bigg\rangle
\\
& + \, \sum_{i=1}^N \bigg\langle \big(\mathrm d \hspace{0.5pt} \J^2_i\big)_{ (\vec u_p, \vec v_p,\xi_p)}\big(\vec z, \vec w,\eta-\xi_p\big),\, \psi_{pi} \bigg\rangle.
\end{split}
\eeq  
\end{proposition}

\bp By Lemmas \ref{pr7}-\ref{le8}, $\I_p$ is Frech\'et differentiable and $\J$ is a continuously Frech\'et differentiable submersion on $\mathscr{X}^p(\Om)$. Also, the set $\mathscr{X}^p_M(\Om)$ is convex and has non-empty interior (with respect to the strong topology). Hence, the assumptions of the generalised Kuhn-Tucker theorem hold true (see e.g.\ \cite[p.\ 417-418, Corollary 48.10 \& Theorem 48B]{Z}). Hence, there exists a Lagrange multiplier
\[
\begin{split}
\Lambda_p \ \in \ & \Big(\big({\mathrm{W}}^{1,\frac{m}{m-2}} (\Om;\R^{2\by N})\big)^* \by \big({\mathrm{W}}^{1,\frac{p}{p-1}} (\Om;\R^{2\by N})\big)^*\Big)^*
\end{split}
\]
which can be identified with a pair of functions
\[
\big(\vec \phi_p, \vec \psi_p \big) \ \in {\mathrm{W}}^{1,\frac{m}{m-2}} (\Om;\R^{2\by N}) \by {\mathrm{W}}^{1,\frac{p}{p-1}} (\Om;\R^{2\by N})
\]
such that, $\big(\vec u_p, \vec v_p, \xi_p\big)$ satisfies
\beq \label{5.19}
\begin{split}
  \frac{1}{p}  \big(\mathrm d \hspace{0.5pt} \I_p   &  \big)_{ (\vec u_p, \vec v_p,\xi_p)}    \Big(\vec z -\vec u_p, \, \vec w -\vec v_p,\, \eta-\xi_p\Big)
 \\
 & \geq \, \sum_{i=1}^N \bigg\langle \big(\mathrm d \hspace{0.5pt} \J^1_i\big)_{ (\vec u_p, \vec v_p,\xi_p)} \Big(\vec z -\vec u_p, \, \vec w -\vec v_p,\, \eta-\xi_p\Big),\, \phi_{pi}   \bigg\rangle
\\
& \ \ \ \, + \, \sum_{i=1}^N \bigg\langle \big(\mathrm d \hspace{0.5pt} \J^2_i\big)_{ (\vec u_p, \vec v_p,\xi_p)} \Big(\vec z -\vec u_p, \, \vec w -\vec v_p,\, \eta-\xi_p\Big),\, \psi_{pi} \bigg\rangle,
\end{split}
\eeq
for any $(\vec z, \vec w,\eta) \in \mathscr{X}^p_M(\Om)$. Since the convex set $\mathscr{X}^p_M(\Om)$ can be written as the cartesian product of the vector spaces 
\[
{\mathrm{W}}^{1,\frac{m}{2}} (\Om;\R^{2\by N}) \by {\mathrm{W}}^{1,p} (\Om;\R^{2\by N})
\]
with the convex set ${\mathrm{L}}^p(\Om,[0,M]) $ (see \eqref{5.10}), by replacing $\vec z$ by $\vec z +\vec u_p$ and $\vec w$ by $\vec w +\vec v_p$ in \eqref{5.19}, we arrive at \eqref{5.17}. The proof of the Proposition is complete. 
\ep

By Proposition \ref{pr9} we deduce that the variational inequality takes the form \eqref{5.18} below, as a direct consequence of \eqref{5.1}, \eqref{5.3}, \eqref{5.2}, \eqref{5.12}-\eqref{5.15}. 

\begin{corollary} \label{cor10}
In the setting of Proposition \ref{pr9}, in view of the form of the Frech\'et derivatives of $\I_p$ and $\J$, the variational inequality \eqref{5.17} takes the form
\beq
\label{5.18}
\begin{split}
& \int_{\p\Om} \vec w : \mathrm d [\vec \nu_p(\vec v_p)]\, + \, \al \int_\Om (\eta-\xi_p)\, \mathrm d [\mu_p(\xi_p)]
\\
& \geq \, \sum_{i=1}^N \bigg\{\int_\Om\Big[ \A :(\D z_i^\top \D \phi_{pi}) \, + (\K z_i)\cdot \phi_{pi} \Big] \, \mathrm{d}\mL^n 
\\
& \ \ \ + \int_{\p\Om}(\ga z_i)\cdot \phi_{pi} \, \mathrm{d}\mH^{n-1}\bigg\}
 + \sum_{i=1}^N \bigg\{\int_\Om \Big[ \B :(\D w_i^\top \D \psi_{pi}) 
\\
&\ \ \ +  \Big(\L w_i- \M \big((\eta-\xi_p) u_{pi} +\xi_p z_i \big)\Big) \cdot \psi_{pi} \Big] \, \mathrm{d}\mL^n 
\, + \int_{\p\Om} (\ga w_i)\cdot \psi_{pi} \, \mathrm{d}\mH^{n-1} \bigg\},
\end{split}
\eeq
for any $ (\vec z, \vec w, \eta  ) \in \mathscr{X}^p_M(\Om)$. 
\end{corollary}

We conclude this section by obtaining the further desired information on the variational inequality \eqref{5.18}.

\begin{lemma} \label{le11}
In the setting of Corollary \ref{cor10}, the variational inequality \eqref{5.18} for the constrained minimiser $\big(\vec u_p, \vec v_p ,\xi_p \big)$ is equivalent to the triplet of relations \eqref{5.20}-\eqref{5.22}.
\end{lemma}

\bp The inequality \eqref{5.20} follows by setting $\vec z =\vec w =0$ in \eqref{5.18}, and recalling the definition of Radon-Nikodym derivative of the absolutely continuous measure $\mu_p(\xi_p)$. The identity \eqref{5.21} follows by setting $\eta =\xi_p$ and $\vec z =0$ in \eqref{5.18} and by recalling that ${\mathrm{W}}^{1,p} (\Om;\R^{2\by N})$ is a vector space, so the inequality we obtain in fact holds for both $\pm w$. Finally, the identity \eqref{5.22} follows by setting $\eta =\xi_p$ and $\vec w =0$ in \eqref{5.18} and by recalling again that ${\mathrm{W}}^{1,\frac{m}{2}} (\Om;\R^{2\by N})$ is a vector space, so the inequality holds for both $\pm z$.
\ep

We conclude by establishing our last main result.

\BPT \ref{th14}. We first show that for any $p>n$ and any 
\[
(\vec v, \xi) \in {\mathrm{W}}^{1,p}(\Om;\R^{2\by N}) \by {\mathrm{L}}^p(\Om), 
\]
we have the next total variations bounds for the measures \eqref{5.3}-\eqref{5.2}:
\begin{align}
\big\| \vec \nu_p(\vec v)\big\|(\p\Om)\, &\leq\, N, \label{5.4}
\\
\big\| \mu_p(\xi)\big\|(\Om)\, & \leq\, 1.  \label{5.5}
\end{align}
To see \eqref{5.4}-\eqref{5.5}, note that for any $i\in\{1,...,N\}$ by H\"older inequality we have
\[
\begin{split}
\big\| \nu_{p_i}(\vec v) \big\|(\p\Om) \, & \leq\, \frac{ \av_{\mathrm{B}_i}\big(|v_i-\tilde v_i|_{(p)}\big)^{p-2} |v_i-\tilde v_i| \, \mathrm d  \mH^{n-1}  }{\left(\, \av_{\mathrm{B}_i}\big(\big| v_i -\tilde v_i \big|_{(p)} \big)^p\, \mathrm d \mH^{n-1}\right)^{\!\!\frac{p-1}{p}}}
\\
 &  \leq\, \frac{ \av_{\mathrm{B}_i}\big(|v_i-\tilde v_i|_{(p)}\big)^{p-1} \, \mathrm d  \mH^{n-1}  }{\left(\, \av_{\mathrm{B}_i}\big(\big| v_i -\tilde v_i \big|_{(p)} \big)^p\, \mathrm d \mH^{n-1}\right)^{\!\!\frac{p-1}{p}}}
\\
&\leq \, 1.
\end{split}
\]
Similarly, we estimate
\[
\begin{split}
\big\| \mu_{p}(\xi) \big\|(\Om) \, & \leq\, \frac{ \av_{\Om}\big(|\xi|_{(p)}\big)^{p-2} |\xi| \, \mathrm d  \mL^{n}  }{\left(\, \av_{\Om}\big(| \xi |_{(p)} \big)^p\, \mathrm d \mL^{n}\right)^{\!\!\frac{p-1}{p}}}
\, \leq\, \frac{ \av_{\Om}\big(|\xi|_{(p)}\big)^{p-1} \, \mathrm d  \mL^{n}  }{\left(\, \av_{\Om}\big(| \xi |_{(p)} \big)^p\, \mathrm d \mL^{n}\right)^{\!\!\frac{p-1}{p}}}
\, \leq \, 1.
\end{split}
\]
Further, by the sequential weak* compactness of the corresponding spaces of Radon measures, the estimates \eqref{5.4}-\eqref{5.5} imply the existence of a subsequence $(p_j)_1^\infty$ and of the claimed limiting measures $(\mu_\infty,\vec \nu_\infty)$ in \eqref{6.2}. Note that the non-negativity of $\mu_\infty$ follows from that of $\xi_p$ and hence of $\mu_p(\xi_p)$.

Next, we prove for later use the estimate
\beq
\label{6.8}
\liminf_{p_j \to \infty}\int_\Om \xi_p \, \mathrm d[\mu_p(\xi_p)] \,\geq\,  \|\xi_\infty\|_{{\mathrm{L}}^\infty(\Om)}.
\eeq
To see \eqref{6.8}, we argue as follows. First, note that if $\xi_\infty = 0$ a.e.\ on $\Om$, then by the positivity of $\mu_p$ and $\xi_p$ we trivially have  
\[
\liminf_{p_j \to \infty}\int_\Om \xi_p \, \mathrm d[\mu_p(\xi_p)] \, \geq\,  0\, =\, \|\xi_\infty\|_{{\mathrm{L}}^\infty(\Om)}
\]
and hence \eqref{6.8} ensues. Therefore, we may assume $\|\xi_\infty\|_{{\mathrm{L}}^\infty(\Om)}>0$. Next, note that by \eqref{4.7} we have
\[
\begin{split}
\int_\Om \xi_p \, \mathrm d[\mu_p(\xi_p)] \,&=\, \av_\Om \frac{\big(|\xi_p|_{(p)}\big)^{p-2}\, |\xi_p|^2}{ \big(\|\xi_p\|_{\dot{\mathrm{L}}^p(\Om)}\big)^{p-1}}\,  \mathrm d \mL^n
\\
&=\, \av_\Om \frac{|\xi_p|_{(p)}^{p}}{\big(\|\xi_p\|_{\dot{\mathrm{L}}^p(\Om)}\big)^{p-1}}\,  \mathrm d \mL^n\, - \,
\frac{1}{p^2}\av_\Om \frac{\big(|\xi_p|_{(p)}\big)^{p-2}}{\big(\|\xi_p\|_{\dot{\mathrm{L}}^p(\Om)}\big)^{p-1}}\,  \mathrm d \mL^n,
\end{split}
\]
which by H\"older inequality gives
\[
\begin{split}
\int_\Om \xi_p \, \mathrm d[\mu_p(\xi_p)] \,&=\, \|\xi_p\|_{\dot{\mathrm{L}}^p(\Om)}  - \,
\frac{1}{p^2}\big(\|\xi_p\|_{\dot{\mathrm{L}}^p(\Om)}\big)^{1-p} \av_\Om \big(|\xi_p|_{(p)}\big)^{p-2}\,  \mathrm d \mL^n
\\
&\geq \, \|\xi_p\|_{\dot{\mathrm{L}}^p(\Om)}  - \,
\frac{1}{p^2 \, \|\xi_p\|_{\dot{\mathrm{L}}^p(\Om)}}.
\end{split}
\]
Hence, for any $k \geq 1$ fixed and $p\geq k$, we have
\[
\begin{split}
\int_\Om \xi_p \, \mathrm d[\mu_p(\xi_p)] \,\geq \, \|\xi_p\|_{\dot{\mathrm{L}}^k(\Om)}  - \,
\frac{1}{p^2 \, \|\xi_p\|_{\dot{\mathrm{L}}^k(\Om)}}.
\end{split}
\]
Since by Theorem \ref{th4} we have $\xi_p \weak \xi_\infty$ in ${\mathrm{L}}^k(\Om)$ for any $k\in (1,\infty)$, by the weak lower semi-continuity of the convex functional $\|\cdot \|_{\dot{\mathrm{L}}^k(\Om)}$ on ${\mathrm{L}}^k(\Om)$, it follows that
\[
\begin{split}
\liminf_{p_j \to \infty}\int_\Om \xi_p \, \mathrm d[\mu_p(\xi_p)] \, &\geq \, \liminf_{p_j \to \infty} \, \|\xi_p\|_{\dot{\mathrm{L}}^k(\Om)} - \, \bigg(\underset{p_j \to \infty}{\limsup}\, \frac{1}{p^2}\bigg) \frac{ 1 }{\underset{p_j \to \infty}{\liminf} \, \|\xi_p\|_{\dot{\mathrm{L}}^k(\Om)}}
\\
&\geq \, \|\xi_\infty\|_{\dot{\mathrm{L}}^k(\Om)}.
\end{split}
\]
We therefore discover \eqref{6.8} by letting $k\to \infty$.

\smallskip

Now we proceed with establishing (I) and (II) of the theorem.

\smallskip

\noi (I) Suppose that $C_\infty=0$. Then, we have 
\beq
\label{6.9}
\ \ \  \big( \vec \phi_p, \vec \psi_p\big) \larrow \big(\vec 0, \vec 0 \big)\ \ \text{ in }  {\mathrm{W}}^{1,\frac{m}{m-2}} (\Om;\R^{2\by N})\! \by \mathrm{BV} (\Om;\R^{2\by N})
\eeq
as $p_j \to \infty$, where $\big(\vec \phi_p, \vec \psi_p \big)$ are the Lagrange multipliers associated with the constrained minimisation problem \eqref{4.15B}. In view of \eqref{6.8} and \eqref{4.16}, the inequality \eqref{5.20} implies
\beq
\label{6.10}
\al \! \int_\Om \eta \, \mathrm d [\mu_p(\xi_p)] \, + \, \sum_{i=1}^N \int_\Om (\eta-\xi_p)\big(\M  u_{pi}  \big)  \cdot \psi_{pi} \, \mathrm{d}\mL^n \, \geq\, o(1)_{p_j\to \infty} +\, \al\|\xi_\infty\|_{{\mathrm{L}}^\infty(\Om)}.
\eeq
for any $\eta \in {\mathrm{C}}^0_0(\Om,[0,M])$. Note now that H\"older's inequality gives
\[
\int_\Om \Big|(\eta-\xi_p) \M  u_{pi} \Big|^{\frac{m}{2}} \, \mathrm{d}\mL^n \, \leq\, \bigg(\int_\Om \big| \eta-\xi_p\big|^m \, \mathrm{d}\mL^n \bigg)^{\frac{1}{2}} \bigg(\int_\Om \big| \M  u_{pi}\big|^m \, \mathrm{d}\mL^n \bigg)^{\frac{1}{2}}
\]
and by \eqref{3.4} and \eqref{4.16} the right hand side of the above estimate is bounded uniformly in $p$.  Hence, by \eqref{6.2}, \eqref{6.9} and the weak-strong continuity of the duality pairing between ${\mathrm{L}}^{\frac{m}{2}}(\Om)$ and ${\mathrm{L}}^{\frac{m}{m-2}}(\Om)$, by letting $p\to \infty$ along the sequence $(p_j)_1^\infty$, \eqref{6.10} yields
\[
\al \! \int_\Om \eta \, \mathrm d \mu_\infty \, \geq\,  \al\|\xi_\infty\|_{{\mathrm{L}}^\infty(\Om)},
\]
for any $\eta \in {\mathrm{C}}^0_0(\Om,[0,M])$. Hence, if $\al>0$, we see that $\xi_\infty=0$ a.e.\ on $\Om$. Again by \eqref{6.2} and \eqref{6.9}, by passing to the limit as $p_j \to \infty$ in \eqref{5.21}, we obtain
\[
\begin{split}
 \int_{\p\Om} \vec w : \mathrm d \vec \nu_\infty \, =\,0 \, =\, \lim_{p_j\to \infty}\sum_{i=1}^N & \bigg\{  \int_\Om \Big[ \B :(\D w_i^\top \D \psi_{pi}) \, +  \big(\L w_i\big) \cdot \psi_{pi} \Big] \, \mathrm{d}\mL^n 
\\
& \ + \int_{\p\Om} (\ga w_i)\cdot \psi_{pi} \, \mathrm{d}\mH^{n-1} \bigg\},
\end{split}
\]
for any $\vec w  \in {\mathrm{C}}^1_0 (\overline{\Om};\R^{2\by N})$. Therefore, $\vec \nu_\infty =\vec 0$, as claimed.

\ms

\noi (II) Suppose now that $C_\infty>0$. Then, the desired relations \eqref{6.3}-\eqref{6.5} would follow directly from \eqref{6.10} and \eqref{5.21}-\eqref{5.22} by rescaling $\big(\vec \phi_p ,\vec \psi_p \big)$ and passing to the limit as $p_j \to \infty$ since the rescaled multipliers $\big(\vec \phi_p/C_p ,\vec \psi_p/C_p \big)$ are bounded in the product space
\[
{\mathrm{W}}^{1,\frac{m}{m-2}} (\Om;\R^{2\by N}) \by \mathrm{BV} (\Om;\R^{2\by N})
\]
and therefore the sequence is sequentially weakly* compact, once we justify the convergence
\beq
\label{6.11}
\int_\Om \xi_p \big( \M u_{pi}\big)\cdot \frac{\psi_{pi}}{C_p}\, \mathrm d \mL^n \larrow \int_\Om \xi_\infty \big( \M u_{\infty i}\big)\cdot \psi_{\infty i}\, \mathrm d \mL^n,
\eeq
as $p_j \to \infty$. To this end, we estimate
\beq
\label{6.12}
\begin{split}
\bigg| \int_\Om  \xi_p \big( \M u_{pi}\big)\cdot \frac{\psi_{pi}}{C_p} &\, \mathrm d \mL^n - \int_\Om \xi_\infty \big( \M u_{\infty i}\big)\cdot \psi_{\infty i}\, \mathrm d \mL^n \bigg|
\\
& \leq\, \int_\Om  |\xi_p|\bigg| \big( \M u_{pi}\big)\cdot \frac{\psi_{pi}}{C_p}\, - \big( \M u_{\infty i}\big)\cdot \psi_{\infty i}\bigg|\, \mathrm d \mL^n 
\\
&\ \ \ \ + \bigg| \int_\Om (\xi_p-\xi_\infty) \big[\big( \M u_{\infty i}\big)\cdot \psi_{\infty i}\big]\, \mathrm d \mL^n \bigg| .
\end{split}
\eeq
Note now that by Theorem \ref{th4} we have $\xi_p \weak \xi_\infty$ in ${\mathrm{L}}^q(\Om)$ for any $q\in (1,\infty)$ as $p_j\to \infty$. In view of \eqref{6.12}, to conclude with \eqref{6.11} we need to show that there exists $t>1$ such that
\begin{align}
&\Big\|\big( \M u_{\infty i}\big)\cdot \psi_{\infty i}\Big\|_{{\mathrm{L}}^t(\Om)} <\, \infty,
\label{6.14}
\\
\label{6.13}
\lim_{p_j\to\infty} & \bigg\| \big( \M u_{pi}\big)\cdot \frac{\psi_{pi}}{C_p}\, - \big( \M u_{\infty i}\big)\cdot \psi_{\infty i}\bigg\|_{{\mathrm{L}}^t(\Om)}=\, 0.
\end{align}
To this end, recall that by Theorem \ref{th4} we have $u_{pi}\weak u_{\infty i}$ in ${\mathrm{W}}^{1,\frac{m}{2}}(\Om;\R^2)$ as $p_j\to \infty$. Thus, by standard compactness arguments we have $u_{pi}\larrow u_{\infty i}$ in ${\mathrm{L}}^{\frac{m}{2}}(\Om;\R^2)$ and hence a.e.\ on $\Om$ as $p_j\to \infty$ (perhaps up to a further subsequence). 
Without loss of generality we suppose that $m<2n$. (If $m\geq 2n$, replace $m$ by any number $\tilde m <2n$.) Since
\[
\Big(\frac{m}{2}\Big)^* =\, \frac{nm}{2n-m}
\]
by the Sobolev inequalities, we have that $(u_{pi})_p$ is bounded in the space ${\mathrm{L}}^{\frac{nm}{2n-m}}(\Om;\R^2)$. By this bound and the a.e.\ convergence  as $p_j\to \infty$, the Vitali convergence theorem (e.g.\ \cite{FL}) implies
\beq
\label{6.15}
\text{$u_{pi}\larrow u_{\infty i}$ \ \ in ${\mathrm{L}}^{\frac{nm(1-\e)}{2n-m}}(\Om;\R^2)$ as $p_j\to \infty$, \ for any }\e \in(0,1).
\eeq
We now argue similarly for $(\psi_{pi}/C_p)_{p_j=1}^\infty$. Since the sequence is bounded and weakly* convergent in $\mathrm{BV}(\Om;\R^2)$, by the compactness of the embedding of $\mathrm{BV}$ in ${\mathrm{L}}^1$ we have $\psi_{pi}/C_p\larrow \psi_{\infty i}$ in ${\mathrm{L}}^{1}(\Om;\R^2)$ and hence a.e.\ on $\Om$ as $p_j\to \infty$ (perhaps up to a further subsequence). By the Sobolev inequality, we have that $(\psi_{pi}/C_p)_p$ is bounded in the space ${\mathrm{L}}^{\frac{n}{n-1}}(\Om;\R^2)$. By this bound and the a.e.\ convergence  as $p_j\to \infty$, the Vitali convergence theorem implies
\beq
\label{6.16}
\text{$\frac{\psi_{pi}}{C_p}\larrow \psi_{\infty i}$ \ \ in ${\mathrm{L}}^{\frac{n(1-\e)}{n-1}}(\Om;\R^2)$ as $p_j\to \infty$, \ for any }\e \in(0,1).
\eeq
Fix now $r,s\geq 1$. By H\"older's inequality we have
\beq
\label{6.17}
\int_\Om \left| \big( \M u_{pi}\big)\cdot \frac{\psi_{pi}}{C_p} \right|^r \mathrm d \mL^n \, \leq\, \left(\int_\Om  \big| \M u_{pi}\big|^{rs}\, \mathrm d \mL^n \right)^{\!1/s} \left(\int_\Om  \Big| \frac{\psi_{pi}}{C_p}\Big|^{rs'}\, \mathrm d \mL^n \right)^{\!1/s'}.
\eeq
We now select
\[
r\,:=\, \frac{1}{1-\e}
\]
and
\[
\e \,:=\, 1 - \sqrt \kappa
\]
where
\[
 \kappa\,:=\, \max\left\{ 1-2\frac{m-n}{nm} \, ,\, \frac{2n-m}{nm} \right\}.
\]
In order to show the above choices are admissible, we need to prove that $0<\kappa<1$. To this aim, we readily have that $\kappa>0$ because $m<2n$. Next, note we also have
\[
1-2\frac{m-n}{nm} \,<\,1,
\]
because $m>n$. Further, we have the equivalence
\[
\frac{2n-m}{nm}\,<\, 1 \ \ \Longleftrightarrow \ \ m\,>\, \frac{2n}{n+1}
\]
and given that $m>n$ and $\frac{2n}{n+1}<2$, we deduce that $\kappa<1$. Next, we choose
\[
s\,:=\, \left(\frac{nm}{2n-m}\right)\left(\frac{1-\e}{r}\right).
\]
Then, our earlier choices of $r,\e,\ka$ imply that $s$ is an admissible choice, since
\[
s\, =\, \frac{nm \kappa}{2n-m}\, \geq\, \left(\frac{nm }{2n-m}\right)\left(\frac{2n-m}{nm}\right) = \,1.
\]
Note further that
\[
rs' \,\leq\, \frac{n(1-\e)}{n-1}
\]
as the above inequality can be easily seen that is equivalent to
\[
 \kappa\,\geq\, 1-2\frac{m-n}{nm}
\]
and the latter inequality is true by the definition of $\kappa$. In conclusion, by H\"older's inequality and the above arguments, \eqref{6.17} yields
\beq
\label{6.18}
\int_\Om \left| \big( \M u_{pi}\big)\cdot \frac{\psi_{pi}}{C_p} \right|^r  \mathrm d \mL^n  \leq \left(\, \int_\Om  \big| \M u_{pi}\big|^{ \frac{nm(1-\e)}{2n-m} }\, \mathrm d \mL^n \! \right)^{\!\!\frac{1}{s}}  \left( \int_\Om  \Big| \frac{\psi_{pi}}{C_p}\Big|^{\frac{n(1-\e)}{n-1}}\, \mathrm d \mL^n \!\right)^{\!\!\frac{r(n-1)}{n(1-\e)}}.
\eeq
In view of \eqref{6.15}-\eqref{6.16}, \eqref{6.14} ensues from \eqref{6.18} for any $t \in (1,r)$. Finally, \eqref{6.13} also follows from  \eqref{6.18} and the Vitali convergence theorem, as from \eqref{6.15}-\eqref{6.16} we already know
\[
\big( \M u_{pi}\big)\cdot \frac{\psi_{pi}}{C_p}\larrow \big( \M u_{\infty i}\big)\cdot \psi_{\infty i} \ \ \text{ a.e. on }\Om,
\]
as $p_j\to \infty$, because $\M \in {\mathrm{L}}^\infty(\Om;\R^{2\by 2})$. The theorem ensues.
\qed

\ms

\bibliographystyle{amsplain}

\end{document}